\documentclass[10pt,b5paper]{article}
\usepackage[T1]{fontenc}
\usepackage[utf8]{inputenc}
\usepackage[english]{babel}
\usepackage{amsmath,amssymb,amsthm}
\usepackage{lmodern}
\usepackage{geometry}
\usepackage{url}
\usepackage{titlesec}
\usepackage{enumitem}
\titleformat{\paragraph}[runin]{\normalfont\itshape}{}{0pt}{}[.]
\titlespacing*{\paragraph}{0pt}{0.7em}{0.5em}

\newtheorem{theorem}{Theorem}[section]
\newtheorem{lemma}[theorem]{Lemma}
\newtheorem{claim}{Claim} 
\newtheorem{proposition}[theorem]{Proposition}
\newtheorem{corollary}[theorem]{Corollary}
\theoremstyle{definition}
\newtheorem{definition}[theorem]{Definition}
\theoremstyle{remark}
\newtheorem{remark}[theorem]{Remark}
\numberwithin{equation}{section}

\DeclareMathOperator{\E}{E}
\DeclareMathOperator{\Var}{Var}
\DeclareMathOperator{\Cov}{Cov}
\DeclareMathOperator{\Corr}{Corr}
\DeclareMathOperator{\Bin}{Bin}
\DeclareMathOperator{\Mult}{Mult}
\newcommand{\Z}{\mathbb{Z}}
\newcommand{\R}{\mathbb{R}}
\newcommand{\hadj}{h_{\mathrm{adj}}}

\title{Discrete total variation under the multinomial,\\Part 4:\\ 
variance under a general distribution}
\author{
Nikola Bani\'c\\
Gideon Brothers, 10000 Zagreb, Croatia\\
\texttt{nbanic@gmail.com}
\and
Neven Elezovi\'c\\
Department of Applied Mathematics,\\
Faculty of Electrical Engineering and Computing,\\
University of Zagreb, 10000 Zagreb, Croatia\\
\texttt{neven.elezovic@fer.hr}
}
\date{\today}

\begin{document}

\maketitle

\begin{abstract}
We study the variance of the reference-centered discrete total variation
$d_{\mathbf p}=\sum_{j=1}^{n-1}|X_{j+1}-X_j-N(p_{j+1}-p_j)|$ of $X\sim\mathrm{Mult}(N;\mathbf p)$ for fixed $n\ge3$
and fixed $\mathbf p$ with all $p_i>0$, as $N\to\infty$. We prove $\mathrm{Var}(d_{\mathbf p})=N\,v(\mathbf p)+O(1)$,
where $v(\mathbf p)$, the variance of the sum of absolute values of the limiting Gaussian edge vector, is expressed
through the exact multinomial edge variances and correlations and Nabeya's covariance formula for folded bivariate
normal variables. Unlike in the uniform case, non-adjacent edge pairs generally contribute to $v(\mathbf p)$. The
$O(1)$ remainder is controlled by a lattice Edgeworth expansion and Poisson summation; it oscillates with the lattice
phases $\{N(p_{j+1}-p_j)\}$, and for $n=3,4$, where a pair of edges exhausts the bins, the summation runs over a
sublattice coset. The formula reduces to the known uniform result (the classical Hoel--Kamat variance), gives a closed
form when the reference law is a symmetric comb, and has a shape-independent large-bin limit for positive $C^1$ profiles.
For the statistic $d_{\mathbf r}$ centered at a reference law $\mathbf r\ne\mathbf p$ the absolute values linearize: when
every edge is mismatched the variance is that of a linear form in the counts up to an exponentially small error, and in
mixed (triangular-array) regimes it is given to $o(N)$ by bivariate folded-normal covariances. The two-bin formulas are
recorded from the companion preprints, and exact finite-$N$ reductions for $n\le4$ are used for validation.
\end{abstract}



\section{Introduction}

	This paper computes the \emph{variance} of the reference-centered discrete total variation $d_{\mathbf p}$ in
	the matched case $X\sim\Mult(N;\mathbf p)$. Three earlier installments treat the expectation and variance under
	the uniform null~\cite{banic_elezovic_p1,banic_elezovic_p2} and the expectation under a general
	law~\cite{banic_elezovic_p3}; a final part treats the cyclic statistic $d^{c}$~\cite{banic_elezovic_p5}. Here the
	discrete total variation of a histogram is the $\ell_1$ norm of its discrete derivative,
	$\sum_{i=1}^{n-1}|x_{i+1}-x_i|$ --- the unnormalized mean absolute successive difference~\cite{kamat1953}, the $\ell_1$
	analogue of von Neumann's mean-square successive difference~\cite{vonneumann1941} (the Durbin--Watson
	numerator~\cite{durbin_watson1950}) --- and not the total-variation distance between probability measures. With
	the mean of $d_{\mathbf p}$ in hand from Part~3~\cite[Theorem~4.3]{banic_elezovic_p3}, the variance is what remains
	before the statistic can be calibrated against a non-uniform reference without simulation --- the setting posed as
	future work in~\cite{banic_elezovic_ct} (a non-uniform null via multinomial weights) and
	in~\cite{banic_elezovic_tvor} (an analysis of variance of the discrete total variation). Centering at the null is
	what keeps the comb test blind to which bins a comb raises: without it, on a non-uniform null every edge with a
	nonzero expected slope locks its sign, and the uncentered variation collapses to a fixed signed sum of counts that
	detects one comb sign and anti-detects the other (\cite[\S6]{banic_elezovic_p3}, paragraph \emph{The uncentered
	statistic}); Theorem~\ref{thm:Var-dXY} below is the variance counterpart of that collapse.

	\medskip\noindent\textbf{Main result.} For fixed $n\ge3$ and fixed $\mathbf p$ in the open simplex,
	\[
		\Var(d_{\mathbf p})=\Bigl(1-\frac2\pi\Bigr)\sum_{j=1}^{n-1}\Var(D_j)
		+\frac4\pi\sum_{1\le j<k\le n-1}\sigma_j\sigma_k\,\Phi(\rho_{j,k})+O(1),
	\]
	with the exact multinomial edge variances $\Var(D_j)=\sigma_j^2$ and correlations $\rho_{j,k}$ of
	Lemma~\ref{lem:cov-Dj-multinomial} and Nabeya's folded covariance factor
	$\Phi(\rho)=\sqrt{1-\rho^2}+\rho\arcsin\rho-1$ (Theorem~\ref{thm:Var-TVC} for $n\ge5$,
	Corollary~\ref{cor:n34} for $n=3,4$; the two-bin case is sharper and separate, \S\ref{ssec:small-n2}). The leading
	coefficient is the variance of the limit law $\sum_j|Y_j|$ of $d_{\mathbf p}/\sqrt N$ (Corollary~\ref{cor:weak-limit});
	the content of the theorem is the $O(1)$ remainder, its lattice oscillation, and the location of the small-$n$
	coset oscillations at order $N^{-1}$. The remainder is a fixed-$n$, fixed-$\mathbf p$ statement: all large-$n$
	statements refer to the iterated limit $N\to\infty$ then $n\to\infty$ of the leading coefficient
	(Remark~\ref{rem:VarTVC-scope}); the pair-boundary configurations excluded by $q_0>0$ at $n=3,4$ are discussed in
	Remark~\ref{rem:pair-boundary}.

	Under a general matched law the edge variables $\tilde D_j$ have unequal variances and
	correlations $\rho_{j,k}$ of either sign, and the variance of $d_{\mathbf p}$ is assembled from
	the covariances of \emph{folded} normals -- the absolute values $|\tilde D_j|$. The key
	analytic tool is Nabeya's identity, which gives the covariance of two
	folded normals in closed form as a function of their correlation; combined with the
	exact multinomial edge correlations it yields the leading variance. The paper then treats
	three consequences of this formula. The symmetric \emph{comb}, the period-$2$ pattern targeted by
	the comb test~\cite{banic_elezovic_ct} (the discrete total variation itself was introduced as a
	smoothness score in~\cite{banic_elezovic_tvor}), used here as the reference law $\mathbf r=\mathbf p$,
	admits a closed-form null variance. For positive $C^1$ profiles the leading coefficient has a
	shape-independent limit as $n\to\infty$ (Proposition~\ref{prop:smooth-universal});
	\S\ref{ssec:test-calibration-nonuniform} discusses calibration.
	The reference-centered statistic $d_{\mathbf r}$ has a full folded-normal variance, including
	the bivariate folded covariance.

	\S\ref{sec:multinomial-basics} recalls the exact edge covariances; \S\ref{sec:Var-general} proves the main theorem
	through the bivariate $O(1)$ estimate of Proposition~\ref{prop:pair-O1} and isolates the $\kappa_{2,2}$ part of the
	disjoint remainder; \S\ref{sec:profiles-calibration} checks the uniform limit, treats the comb and smooth profiles,
	and discusses calibration; \S\ref{sec:Var-dXY} computes $\Var(d_{\mathbf r})$ under $\mathbf p\ne\mathbf r$;
	\S\ref{sec:small-n} settles $n\le4$.

	These results belong to a series on the exact and asymptotic moment theory of reference-centered
	discrete total variation under the multinomial model, organized by \emph{which} moment is computed and
	\emph{under which} reference law: the expectation~\cite{banic_elezovic_p1} and
	variance~\cite{banic_elezovic_p2} under the uniform null; the
	expectation~\cite{banic_elezovic_p3} and the variance, treated here, under a
	general distribution; and the cyclic statistic $d^{c}$ in Part~5~\cite{banic_elezovic_p5}.

\section{Multinomial distribution: basics}
	\label{sec:multinomial-basics}

	This section records the multinomial facts used below: marginal distributions,
	the two-bin factorization, and exact covariances of differences between
	neighboring bins. These facts are standard (see, e.g.,
	Johnson--Kotz--Balakrishnan~\cite{johnson_kotz_balakrishnan}), but the proofs
	are included to fix notation and provide labeled references for later use.
	The main reduction is that the difference $D_j=X_{j+1}-X_j$ depends only on a
	two-bin marginal, and that this marginal factors into a binomial total and a
	conditional binomial split.

\subsection{Definition and basic properties}
	\label{ssec:mult-def}

	Let $X = (X_1, X_2, \ldots, X_n)$ be a random vector following the
	multinomial distribution
	\begin{equation}\label{eq:mult}
		X \sim \Mult(N;\, p_1, p_2, \ldots, p_n), \qquad
		\sum_{i=1}^n p_i = 1,\quad p_i \ge 0.
	\end{equation}
	We interpret the components $X_i$ as the number of observations falling into the $i$-th bin
	in $N$ independent trials with probabilities $p_i$. The associated
	probability mass function is
	\begin{equation}\label{eq:mult-pmf}
		\Pr(X_1 = k_1, \ldots, X_n = k_n) =
		\frac{N!}{k_1!\,k_2!\cdots k_n!}\,p_1^{k_1}p_2^{k_2}\cdots p_n^{k_n},
		\quad \sum_i k_i = N,
	\end{equation}
	and $0$ otherwise.

	\paragraph{Marginal distribution of a single bin} The marginal distribution of a
	single bin is binomial,
	\begin{equation}\label{eq:mult-marg-1}
		X_i \sim \Bin(N, p_i),
	\end{equation}
	so $\E[X_i] = N p_i$ and $\Var(X_i) = N p_i(1 - p_i)$.

	\paragraph{Covariance structure} The standard computation (via the representation
	$X_i = \sum_{k=1}^N \mathbf 1[Z_k = i]$ with $Z_k$ iid categorical on
	$\{1, \ldots, n\}$ with probabilities $p_i$) gives
	\begin{equation}\label{eq:mult-cov}
		\Cov(X_i, X_j) =
	\begin{cases}
		N p_i (1 - p_i), & i = j, \\
		- N p_i p_j, & i \ne j.
	\end{cases}
	\end{equation}

\subsection{Marginal distributions of subsets}
	\label{ssec:mult-marginals}

	The multinomial distribution has the following fundamental property: the marginal
	of any subset of variables is also multinomial, where the
	``remainder'' is treated as one additional bin.

\begin{proposition}[$k$-bin marginal]\label{prop:k-marginal}
	Let $S \subseteq \{1, \ldots, n\}$, $|S| = k$, and let
	$X_S^c := N - \sum_{i \in S} X_i$ be the number of observations outside $S$.
	The \emph{augmented} vector $\bigl((X_i)_{i \in S},\,X_S^c\bigr)$ under (\ref{eq:mult}) is
	\begin{equation}\label{eq:k-marginal}
		\bigl((X_i)_{i \in S},\,X_S^c\bigr) \sim \Mult\bigl(N;\, (p_i)_{i \in S},\; 1 - \textstyle\sum_{i\in S} p_i\bigr),
	\end{equation}
	a genuine $(k{+}1)$-bin multinomial whose last bin is the ``combined remainder'' $X_S^c$,
	with probability $1 - \sum_{i\in S} p_i$; the marginal $(X_i)_{i\in S}$ alone is its first $k$ coordinates.
\end{proposition}

\begin{proof}
	Writing the corresponding multinomial probability in explicit form
	(\ref{eq:mult-pmf}) and summing over all $(k_i)_{i \notin S}$ with
	fixed $k_i$ for $i \in S$ (under the constraint $\sum_i k_i = N$):
	\[
		\Pr\bigl((X_i)_{i \in S} = (k_i)_{i \in S}\bigr)
		= \frac{N!}{\prod_{i\in S} k_i!}\,\prod_{i\in S}p_i^{k_i}
		\cdot \!\!\sum_{\substack{(k_i)_{i \notin S} \\ \sum_{i\notin S} k_i = m}}\!\!
		\frac{1}{\prod_{i\notin S} k_i!}\prod_{i\notin S} p_i^{k_i},
	\]
	where $m = N - \sum_{i\in S} k_i$. By the multinomial theorem the sum reduces to
	\[
		(\sum_{i\notin S} p_i)^m/m! = (1-\sum_{i\in S}p_i)^m/m!,
	\]
	so
	\[
		\Pr\bigl((X_i)_{i \in S} = (k_i)_{i \in S}\bigr)
		= \frac{N!}{(\prod_{i\in S}k_i!)\,m!}\,
		\Bigl(\prod_{i\in S}p_i^{k_i}\Bigr)\bigl(1-{\textstyle\sum_{i\in S}}p_i\bigr)^m,
	\]
	which is exactly the probability mass function of the multinomial with parameters
	$\bigl(N;\,(p_i)_{i\in S},\,1-\sum_{i\in S}p_i\bigr)$.
\end{proof}

	\paragraph{Special cases} We use the following marginals:
	\begin{itemize}
		\item \textbf{2-bin marginal} $(X_j, X_{j+1})$: $\Mult\bigl(N; p_j, p_{j+1}, 1-p_j-p_{j+1}\bigr)$,
		for $\Var(X_{j+1}-X_j)$ and the factorization below.
		\item \textbf{3-bin marginal} $(X_{j-1}, X_j, X_{j+1})$: for
		$\Cov(D_{j-1},D_j)$, which shares $X_j$.
		\item \textbf{4-bin marginal} $(X_j, X_{j+1}, X_k, X_{k+1})$ for $|k-j| \ge 2$:
		for the disjoint covariance $\Cov(D_j, D_k)$ and (with absolute values)
		$\Cov(|D_j|, |D_k|)$.
	\end{itemize}

\subsection{Factorization of the 2-bin marginal}
	\label{ssec:2bin-factorization}

	For a single edge, the $n$-bin multinomial can be reduced to a two-stage
	binomial experiment: first count how many observations fall in the pair
	$\{j,j+1\}$, then split this total between the two bins. This gives a direct
	binomial representation of $X_{j+1}-X_j$ and is used repeatedly in the
	one-edge calculations.

	Concretely, let
	\begin{equation}\label{eq:Kj-def}
		K_j := X_j + X_{j+1},
	\end{equation}
	the sum of the two bins in the marginal
	$(X_j,X_{j+1})\sim\Mult(N;p_j,p_{j+1},1-p_j-p_{j+1})$.

\begin{proposition}[Factorization of the 2-bin marginal]\label{prop:2bin-factor}
	Assume $p_j+p_{j+1}>0$.
	$K_j \sim \Bin(N, p_j + p_{j+1})$, and conditional on $K_j = k$:
	\begin{equation}\label{eq:Xj1-cond}
		X_{j+1} \mid K_j = k \;\sim\; \Bin(k, \pi_j),
		\qquad \pi_j := \frac{p_{j+1}}{p_j + p_{j+1}}.
	\end{equation}
	Conditional on $K_j=k$, the complement $X_j=k-X_{j+1}$ is
	$\Bin(k,1-\pi_j)$.
\end{proposition}

\begin{proof}
	$K_j$ counts how many of the $N$ observations fall into the union of the bins $\{j, j+1\}$
	(success probability $p_j + p_{j+1}$ per independent trial), so it is
	$\Bin(N, p_j+p_{j+1})$.

	For the conditional distribution, by Proposition~\ref{prop:k-marginal} and (\ref{eq:mult-pmf}):
	\begin{align*}
		\Pr(X_{j+1} = m \mid K_j = k)
		&= \frac{\Pr(X_j = k-m, X_{j+1} = m)}{\Pr(K_j = k)} \\
		&= \frac{\binom{N}{k-m, m, N-k}\,p_j^{k-m}p_{j+1}^m (1-p_j-p_{j+1})^{N-k}}
		{\binom{N}{k}(p_j+p_{j+1})^k (1-p_j-p_{j+1})^{N-k}} \\
		&= \binom{k}{m}\,\pi_j^m (1-\pi_j)^{k-m},
	\end{align*}
	which is the probability mass function of $\Bin(k, \pi_j)$.
\end{proof}

	\paragraph{Implication for the differences $D_j = X_{j+1} - X_j$} Conditional on $K_j = k$:
	$D_j = 2X_{j+1} - k$, so $D_j \mid K_j = k$ has a shifted binomial distribution.
	Under the uniform distribution ($p_j = p_{j+1}$), $\pi_j = 1/2$ so
	$X_{j+1} \mid K_j \sim \Bin(k, 1/2)$, which gives $D_j \mid K_j = k$ as
	the position of a symmetric $\pm 1$ random walk of $k$ steps. This
	identification is central to the derivations that follow.

\subsection{Edge variables and their basic covariances}
	\label{ssec:edge-vars}

	The total variation is built from differences between neighboring bins. We call
	these differences \emph{edge variables}; the individual counts enter mainly
	through their covariances. Before absolute values are applied, their joint
	behavior is given by elementary multinomial covariance identities.

\begin{definition}[Edge variables]\label{def:Dj}
	For a multinomial vector $X = (X_1, \ldots, X_n)$ we introduce the \emph{edge variables}
	\begin{equation}\label{eq:Dj-def}
		D_j := X_{j+1} - X_j, \qquad j = 1, \ldots, n-1,
	\end{equation}
	and their centered versions
	\begin{equation}\label{eq:tilde-Dj-def}
		\tilde D_j := D_j - \E[D_j] = D_j - N(p_{j+1} - p_j).
	\end{equation}
	The matched reference-centered statistic studied in the main part of this paper is
	\begin{equation}\label{eq:TVC-def}
		d_{\mathbf p}:=\sum_{j=1}^{n-1}|\tilde D_j|
		=\sum_{j=1}^{n-1}\bigl|D_j-N(p_{j+1}-p_j)\bigr|.
	\end{equation}
	Equivalently, this is the statistic $d_{\mathbf r}$ with reference law
	$\mathbf r=\mathbf p$. Under the uniform distribution ($p_i \equiv 1/n$),
	$\E[D_j] = 0$, so $\tilde D_j = D_j$ and $d_{\mathbf p}$ reduces to the
	uniform statistic $d$.
\end{definition}

	The covariance structure of the edge variables (\emph{without absolute values})
	is exact and completely determined by (\ref{eq:mult-cov}) and the bilinearity
	of the covariance.

\begin{lemma}[Exact covariances of edge variables]\label{lem:cov-Dj-multinomial}
	\leavevmode\par\noindent
	Under the multinomial $\Mult(N;p_1,\ldots,p_n)$:
	\begin{align}
		\Var(D_j) &= N\bigl[(p_j + p_{j+1}) - (p_j - p_{j+1})^2\bigr], \label{eq:Var-Dj} \\
		\Cov(D_j, D_{j+1}) &= N\bigl[-p_{j+1}(1 + p_j + p_{j+2} - p_{j+1}) + p_j p_{j+2}\bigr], \label{eq:Cov-Dj-adj} \\
		\Cov(D_j, D_k) &= -N(p_{j+1} - p_j)(p_{k+1} - p_k) \quad \text{for } |k - j| \ge 2. \label{eq:Cov-Dj-disj}
	\end{align}
	Centering changes neither the variance nor the covariance: $\Var(\tilde D_j) = \Var(D_j)$,
	$\Cov(\tilde D_j, \tilde D_k) = \Cov(D_j, D_k)$.
\end{lemma}

\noindent
The adjacent form \eqref{eq:Cov-Dj-adj} is symmetric in $p_j$ and $p_{j+2}$, as it
must be under reflection of the bin order. Equivalently,
$\Cov(D_j,D_{j+1})=-N\bigl[p_{j+1}+\delta_j\delta_{j+1}\bigr]$ with $\delta_j=p_{j+1}-p_j$:
the shared-bin term plus the same product that appears for disjoint edges. For general
$\mathbf p$ the uncentered second moment $\E D_j^2=N^2(p_{j+1}-p_j)^2+\Var(D_j)$ already
appears in~\cite[Eq.~(11)]{banic_elezovic_tvor}; the adjacent and disjoint covariances
are new there.

\begin{proof}
	All three claims follow from the bilinearity of the covariance: for random
	variables $A, B, C, D$
	\begin{equation}\label{eq:cov-bilin}
		\Cov(A - B,\, C - D) = \Cov(A, C) - \Cov(A, D) - \Cov(B, C) + \Cov(B, D),
	\end{equation}
	and from the multinomial covariances (\ref{eq:mult-cov}).

	\medskip
	\emph{Proof for $\Var(D_j)$.} Special case of (\ref{eq:cov-bilin}) with
	$A = C = X_{j+1}$, $B = D = X_j$:
	\[
		\Var(D_j) = \Var(X_{j+1}) + \Var(X_j) - 2\Cov(X_j, X_{j+1}).
	\]
	Substituting (\ref{eq:mult-cov}):
	\begin{align*}
		\Var(D_j) &= N p_{j+1}(1-p_{j+1}) + N p_j(1-p_j) + 2N p_j p_{j+1} \\
		&= N\bigl[(p_j+p_{j+1}) - (p_j^2 + p_{j+1}^2 - 2 p_j p_{j+1})\bigr]\\
		&= N\bigl[(p_j+p_{j+1}) - (p_j-p_{j+1})^2\bigr].
	\end{align*}

	\medskip
	\emph{Proof for $\Cov(D_j, D_{j+1})$.} Adjacent edges share the variable
	$X_{j+1}$. We apply (\ref{eq:cov-bilin}) with $A = X_{j+1}$, $B = X_j$,
	$C = X_{j+2}$, $D = X_{j+1}$:
	\[
		\Cov(D_j, D_{j+1}) = \Cov(X_{j+1}, X_{j+2}) - \Var(X_{j+1})
		- \Cov(X_j, X_{j+2}) + \Cov(X_j, X_{j+1}).
	\]
	Substituting (\ref{eq:mult-cov}):
	\begin{align*}
		\Cov(D_j, D_{j+1}) &= -N p_{j+1}p_{j+2} - N p_{j+1}(1-p_{j+1}) + N p_j p_{j+2} - N p_j p_{j+1} \\
		&= N\bigl[-p_{j+1}(1 + p_j + p_{j+2} - p_{j+1}) + p_j p_{j+2}\bigr].
	\end{align*}

	\medskip
	\emph{Proof for $\Cov(D_j, D_k)$, $|k-j| \ge 2$.} The indices $\{j, j+1\}$ and
	$\{k, k+1\}$ are disjoint. Special case of (\ref{eq:cov-bilin}) with
	$A = X_{j+1}$, $B = X_j$, $C = X_{k+1}$, $D = X_k$:
	\begin{align*}
		\Cov(D_j, D_k) &= -N p_{j+1} p_{k+1} + N p_{j+1} p_k + N p_j p_{k+1} - N p_j p_k \\
		&= -N \bigl[p_{j+1}(p_{k+1}-p_k) - p_j(p_{k+1}-p_k)\bigr]\\
		&= -N(p_{j+1}-p_j)(p_{k+1}-p_k).
	\end{align*}

	Centering does not affect covariance:
	\[
		\Cov(\tilde D_j,\tilde D_k)
		=\Cov(D_j-\E D_j,D_k-\E D_k)=\Cov(D_j,D_k).
	\]
\end{proof}

\begin{remark}[Geometric interpretation of the disjoint covariance]
	\label{rem:cov-disj-multinomial}
	Since $\E[D_j] = N(p_{j+1} - p_j)$, formula (\ref{eq:Cov-Dj-disj}) can
	be written more compactly:
	\begin{equation}\label{eq:Cov-Dj-disj-compact}
		\Cov(D_j, D_k) = -\frac{\E[D_j]\,\E[D_k]}{N} \qquad (|k-j| \ge 2).
	\end{equation}
	The relation is a consequence of the deterministic constraint
	$\sum_iX_i=N$, which links all multinomial components. In particular, distant edges
	have nonzero covariance whenever both expected edge differences are nonzero and
	their product does not vanish.

	Under the uniform distribution $p_i\equiv1/n$, one has $\E[D_j]\equiv0$ and
	therefore $\Cov(D_j,D_k)=0$ for $|k-j|\ge2$. The disjoint covariance is zero in
	the uniform case, but in general it is not.
\end{remark}

\begin{corollary}[Uniform case]\label{cor:cov-Dj-uniform}
	Under the uniform multinomial ($p_i \equiv 1/n$):
	\begin{align}
		\Var(D_j) &= \frac{2N}{n}, \label{eq:Var-Dj-uniform} \\
		\Cov(D_j, D_{j+1}) &= -\frac{N}{n}, \label{eq:Cov-adj-uniform} \\
		\Cov(D_j, D_k) &= 0, \quad |k - j| \ge 2. \label{eq:Cov-disj-uniform}
	\end{align}
	The corresponding correlations are $\Corr(D_j, D_{j+1}) = -1/2$,
	$\Corr(D_j, D_k) = 0$ for $|k-j| \ge 2$.
\end{corollary}

\begin{proof}
	Substituting $p_i = 1/n$ into (\ref{eq:Var-Dj})--(\ref{eq:Cov-Dj-disj}):
	$(p_j + p_{j+1}) = 2/n$, $(p_j - p_{j+1})^2 = 0$, $p_{j+1}(1 + p_j + p_{j+2} - p_{j+1}) = (1/n)(1+1/n)$, whose $1/n^2$ cancels against $p_jp_{j+2}=1/n^2$
	and $(p_{j+1}-p_j)(p_{k+1}-p_k) = 0$.
\end{proof}

	\paragraph{Correlations} For adjacent edges we introduce
	\begin{equation}\label{eq:rho-j-mult}
		\rho_j := \frac{\Cov(D_j, D_{j+1})}{\sqrt{\Var(D_j)\Var(D_{j+1})}},
	\end{equation}
	and for a general pair $1 \le j < k \le n-1$
	\begin{equation}\label{eq:rho-jk-mult}
		\rho_{j,k} := \frac{\Cov(D_j, D_k)}{\sqrt{\Var(D_j)\Var(D_k)}}.
	\end{equation}
	For the uniform distribution, $\rho_j \equiv -1/2$ and $\rho_{j,k} \equiv 0$ for
	$|k-j| \ge 2$. For a general non-uniform distribution, both adjacent and distant
	correlations can have either sign. For example, if the middle probability
	$p_{j+1}$ is small while $p_j$ and $p_{j+2}$ are large, the positive term
	$p_jp_{j+2}$ in (\ref{eq:Cov-Dj-adj}) can dominate. For distant pairs, the sign
	is opposite to the sign of $(p_{j+1}-p_j)(p_{k+1}-p_k)$
	(see Remark~\ref{rem:cov-disj-multinomial}).

	Two structural facts deserve emphasis, since the whole variance theory rests on
	them. First, every covariance in Lemma~\ref{lem:cov-Dj-multinomial} is
	\emph{exactly linear in $N$}, so the correlations $\rho_j$, $\rho_{j,k}$ and the
	normalized variances $\Var(D_j)/N$ do not depend on $N$ at all: the leading
	variance coefficient obtained below is a genuine function of $\mathbf p$ alone.
	Second, if $\min_ip_i>0$ then each edge variance is strictly positive
	($\Var(D_j)/N=s_j-\delta_j^2\ge s_j(1-s_j)>0$ for $n\ge3$, where
	$s_j=p_j+p_{j+1}$, $\delta_j=p_{j+1}-p_j$) and every correlation is bounded away
	from $\pm1$: $|\rho_{j,k}|\le1-c(\mathbf p)<1$ for all pairs, because the
	per-trial increment pair is supported on at least three affinely independent
	points with positive probability, so its covariance matrix is nonsingular. In
	particular the Nabeya covariance factor below is always evaluated strictly
	inside its domain.

\subsection{Imported and auxiliary estimates}
	\label{ssec:imported}

The variance theorem and its supporting propositions draw on three asymptotic facts from the
earlier installments of the series and on one elementary Fourier estimate, proved here. We collect them
before they are used, with the hypotheses and conclusions needed below. All are for fixed interior $\mathbf p$ (that is, $p_i\ge p_{\min}>0$ for all
$i$) and $N\to\infty$, with constants depending on $(n,\mathbf p)$ unless stated otherwise.
Throughout we write $s_j:=p_j+p_{j+1}$, $\delta_j:=p_{j+1}-p_j$ and
$\tilde s_j:=\sqrt{s_j-\delta_j^2}=\sqrt{\Var(D_j)/N}$; in the notation of
Part~3~\cite{banic_elezovic_p3} these are $\tau_j$, $\mu_j$ and that paper's $s_j$, respectively.

\begin{lemma}[Supporting estimates]\label{lem:imported}
	\leavevmode
	\begin{enumerate}[label=\textup{(\roman*)},leftmargin=*]
	\item \emph{Per-edge folded-mean expansion} \cite[Lemma~4.1]{banic_elezovic_p3}:
		$\E|\tilde D_j|=\sqrt{2/\pi}\,\sigma_j\bigl(1+O(N^{-1})\bigr)$, where
		$\sigma_j^2=\Var(D_j)$; the $O(N^{-1})$ is uniform for $\mathbf p$ in a compact
		subset of the open simplex.
	\item \emph{Kink second-derivative structure} \cite[Lemma~5.2]{banic_elezovic_p3}: the
		distributional second derivative of $x\mapsto|x-c|\,H(x)$, for $H$ smooth, is
		$|x-c|H''+2\operatorname{sgn}(x-c)H'+2\delta_c H$, the last a Dirac mass at the
		kink $c$.
	\item \emph{Scaled-profile Fourier estimate} (proved below): for $g_N(x)=A_NH_N((x-c)/\sigma)$ with
		$A_N=O(1)$, $H_N$ in a bounded Schwartz family, $\sigma\asymp\sqrt N$, and
		$H_N'(0)=O(N^{-1/2})$,
		\[
			\int_{\R}|x-c|\,g_N(x)\,e^{-i\omega x}\,dx
			=-\frac{2g_N(c)\,e^{-i\omega c}}{\omega^2}
			+O\!\bigl(N^{-1}|\omega|^{-3}\bigr),
		\]
		for $\omega=2\pi r$, $r\ne0$, uniformly in $r$, the constant depending only on the
		Schwartz seminorms of the family. Without the condition on the derivative
		at the kink the next term is only $O(N^{-1/2}|\omega|^{-3})$; in the
		applications below the Edgeworth parity gives $H_N'(0)=O(N^{-1/2})$.
	\item \emph{Uniform small-$n$ boundary benchmarks}
		\cite[Theorem~4.1, Lemma~5.1, Theorem~4.6 and Theorem~5.3]{banic_elezovic_p2},
		for the corresponding uniform configurations at $n=3,4$.
	\end{enumerate}
\end{lemma}

\begin{proof}[Proof of \textup{(iii)}]
	Substituting $x=c+\sigma t$ gives
	\[
		\int_{\R}|x-c|\,g_N(x)e^{-i\omega x}\,dx=A_N\sigma^2e^{-i\omega c}F(\omega\sigma),\qquad
		F(\Omega):=\int_{\R}|t|H_N(t)e^{-i\Omega t}\,dt.
	\]
	Three integrations by parts, using \textup{(ii)} for the
	distributional derivatives of $|t|H_N(t)$ and the jumps $2H_N(0)$ and $4H_N'(0)$ of its first and second
	derivatives at the kink, give
	\[
		F(\Omega)=-2H_N(0)\Omega^{-2}+4iH_N'(0)\Omega^{-3}+O(\Omega^{-4}),
	\]
	the constant depending only on the Schwartz seminorms of the family. Multiplying by $A_N\sigma^2$ and using $g_N(c)=A_NH_N(0)$,
	$\sigma\asymp\sqrt N$, $A_N=O(1)$: the first term is $-2g_N(c)e^{-i\omega c}/\omega^2$, the second is
	$O(H_N'(0)N^{-1/2}|\omega|^{-3})=O(N^{-1}|\omega|^{-3})$ under the derivative condition (and only
	$O(N^{-1/2}|\omega|^{-3})$ without it), and the rest is $O(N^{-1}\omega^{-4})$.
\end{proof}

\noindent
For fixed interior $\mathbf p$ and $n\ge5$ the proof of Proposition~\ref{prop:pair-O1} uses
\textup{(i)--(iii)}; the boundary discussion (Remark~\ref{rem:pair-boundary}) uses
\textup{(iv)}. The estimate \textup{(iii)} is the one whose uniformity class matters most
below; what enters
is a bounded Schwartz family with the kink derivative controlled as above, the
scaling $\sigma\asymp\sqrt N$, and the decay $O(N^{-1}|\omega|^{-3})$
uniform in the nonzero Fourier mode $r$.

\section{General variance formula for $d_{\mathbf p}$}
	\label{sec:Var-general}

	This section extends the variance formula for reference-centered discrete
	total variation from the uniform law to the matched statistic $d_{\mathbf p}$
	under a fixed positive multinomial law. In contrast
	with the uniform case, disjoint edge covariances generally contribute to the
	leading coefficient. Consequently the leading variance is obtained by summing
	the Nabeya folded-normal covariance over all pairs of edges, not only
	over adjacent pairs. Remark~\ref{rem:comb-leading} illustrates the numerical
	size of this disjoint contribution for a comb profile.

	We first assemble the exact edge covariances and Nabeya's identity
	into the leading formula for $\Var(d_{\mathbf p})$. Subsequent sections check the uniform
	limit, treat profile families and calibration, and then turn to the variance of
	the reference-centered statistic $d_{\mathbf r}$ under a sampling law
	$\mathbf p$.

\subsection{Exact covariances and correlations}
	\label{ssec:cov-D-nonuniform}

	The exact second-order structure of $(D_1,\dots,D_{n-1})$ is Lemma~\ref{lem:cov-Dj-multinomial}, with the correlations \eqref{eq:rho-j-mult}--\eqref{eq:rho-jk-mult}; it is linear in $N$ and its correlations are $N$-free and bounded away from $\pm1$ (\S\ref{ssec:edge-vars}).

\subsection{Nabeya identity}
	\label{ssec:sheppard-mardia}

	The covariance of absolute values of correlated centered normal variables is
	given by the Nabeya identity. We use it to convert the limiting
	bivariate normal laws of the edge variables into the leading covariance terms
	in $\Var(d_{\mathbf p})$.

	The identity is proved in the companion paper on the uniform
	variance~\cite[Lemma~3.1]{banic_elezovic_p2} by Price's
	theorem~\cite{price1958} (the non-smooth use of Price's theorem being justified
	by the mollification recorded in \cite[Remark~3.2]{banic_elezovic_p2}), where it
	is stated for general
	variances $\sigma_X,\sigma_Y$ precisely so that it can be reused here: for a
	centered bivariate normal $(X, Y) \sim \mathcal N(\mathbf 0, \Sigma)$ with
	$\sigma_X^2 = \Sigma_{11}$, $\sigma_Y^2 = \Sigma_{22}$, correlation
	$\rho = \Sigma_{12}/(\sigma_X \sigma_Y) \in (-1, 1)$,
	\begin{equation}\label{eq:sheppard-mardia}
		\E[|X| \cdot |Y|] = \frac{2}{\pi}\,\sigma_X \sigma_Y
		\bigl[\sqrt{1 - \rho^2} + \rho \arcsin\rho\bigr]
	\end{equation}
	This is the classical bivariate-normal absolute-moment formula: it goes back
	to Sheppard's arcsine formula~\cite{sheppard1899} and was given in this form
	by Nabeya~\cite{nabeya1951}.
	We write
	\begin{equation}\label{eq:SM-Phi}
		\Phi(\rho):=\sqrt{1-\rho^2}+\rho\arcsin\rho-1,\qquad -1<\rho<1,
	\end{equation}
	for the corresponding Nabeya covariance factor; normal distribution
	functions are denoted by $\Phi_{\mathcal N}$ and $\Phi_{\mathcal N_2}$.

\begin{corollary}\label{cor:cov-abs}
	$\Cov(|X|, |Y|) = \frac{2}{\pi}\sigma_X \sigma_Y\Phi(\rho)$.
\end{corollary}
\begin{proof}
	Each marginal is centered normal, so
	$\E|X|=\sigma_X\sqrt{2/\pi}$ and
	$\E|Y|=\sigma_Y\sqrt{2/\pi}$. Subtracting their product
	$2\sigma_X\sigma_Y/\pi$ from \eqref{eq:sheppard-mardia} gives the result.
\end{proof}

\begin{definition}[Folded-normal profile]\label{def:folded-profile}
	For $Z\sim\mathcal N(0,1)$ and $a\in\R$,
	\[
		G(a):=\E|Z+a|=2\phi(a)+a\bigl(2\Phi_{\mathcal N}(a)-1\bigr),
	\]
	the displayed identity following by splitting the integral at $z=-a$ and
	using $\phi'(z)=-z\phi(z)$,
	the folded-normal mean profile ($\phi,\Phi_{\mathcal N}$ the standard normal
	density and CDF; Leone et al.~\cite{leone1961}); equivalently
	$\E|W|=\sigma\,G(\mu/\sigma)$ for $W\sim\mathcal N(\mu,\sigma^2)$. Its further
	properties are collected in \cite[Definition~6.3, Lemma~6.4]{banic_elezovic_p3}.
\end{definition}

\subsection{Main theorem: asymptotic $\Var(d_{\mathbf p})$}
	\label{ssec:Var-TVC-main}

\begin{theorem}[Asymptotic $\Var(d_{\mathbf p})$ with covariances]
	\label{thm:Var-TVC}
	Fix $\mathbf p$ with $p_j>0$ for all $j$, and $n\ge5$. Then, as $N\to\infty$,
	\begin{align}\label{eq:Var-TVC}
		\Var(d_{\mathbf p})
		&= \Bigl(1 - \frac{2}{\pi}\Bigr)\sum_{j=1}^{n-1}\Var(D_j) \notag\\
		&\quad + \frac{4}{\pi}\sum_{1 \le j < k \le n-1}
		\sigma_j \sigma_k \,\Phi(\rho_{j,k}) + R_N,
	\end{align}
	where $\sigma_j := \sqrt{\Var(D_j)}$, $\Phi$ is defined in \eqref{eq:SM-Phi},
	$\rho_{j,k}$ is given by formula \eqref{eq:rho-jk-mult}, and
	$R_N=O(1)$ as $N\to\infty$, with a constant depending on $(n,\mathbf p)$.
\end{theorem}

Equivalently, $\Var(d_{\mathbf p})=N\,v(\mathbf p)+O(1)$, where
$v(\mathbf p)$ is the displayed leading coefficient divided by $N$. This is an
$N$-free function of $\mathbf p$, by the structural facts of
\S\ref{ssec:edge-vars}.

\begin{remark}[Extension to small numbers of bins]
	The statement is formulated here for $n\ge5$, where every pair of edges leaves at
	least one bin outside the pair and the full-lattice argument of
	Proposition~\ref{prop:pair-O1} applies directly. The same variance formula remains
	valid for all $n\ge3$; the exceptional boundary configurations for $n=3,4$ are
	treated separately in \S\ref{sec:small-n} and the final all-$n$ form is recorded in
	Corollary~\ref{cor:n34}.
\end{remark}

	The remainder is genuinely of order one and in general has no limit: through the lattice
	phases $\{N(p_{j+1}-p_j)\}$ of Claim~\ref{cl:axis} and of \cite[Lemma~4.1]{banic_elezovic_p3}
	it is asymptotically periodic in $N$, with period the common denominator of a rational
	$\mathbf p$; for $\mathbf p=(1,3,2,1,3,2)/12$, for example, $\Var(d_{\mathbf p})-Nv(\mathbf p)$
	is about $-0.14$ for $12\mid N$ and about $-1.05$ for $N\equiv4,8\pmod{12}$ (exact
	computation). Only at the uniform law does it converge (to $-\beta(n)$ of Part~2).

\begin{corollary}[The limit law whose variance this is]\label{cor:weak-limit}
	Fix $\mathbf p$ with $p_j>0$ for all $j$, and let $\mathbf Y=(Y_1,\dots,Y_{n-1})$ be
	the centered Gaussian vector with $\Var(Y_j)=\lim_{N\to\infty}\Var(D_j)/N$ and
	$\Corr(Y_j,Y_k)=\rho_{j,k}$, the correlations of \eqref{eq:rho-jk-mult}.
	Then
	\begin{equation}\label{eq:weak-limit}
		\frac{d_{\mathbf p}}{\sqrt N}\ \overset{d}{\longrightarrow}\ W:=\sum_{j=1}^{n-1}|Y_j|,
		\qquad N\to\infty,
	\end{equation}
	and $\Var(d_{\mathbf p})/N\to\Var(W)$. The leading coefficient of
	Theorem~\ref{thm:Var-TVC} is exactly $\Var(W)$: the variance computed there is the
	variance of the identified limit law $W$, the sum of the absolute values of a
	correlated centered Gaussian vector.
\end{corollary}

\begin{proof}
	The centered vector $(\tilde D_1,\dots,\tilde D_{n-1})/\sqrt N$ is a fixed linear
	image of the multinomial count vector, hence asymptotically centered Gaussian with
	the stated variances and correlations by the multivariate central limit theorem
	(centering leaves the covariance unchanged, \S\ref{ssec:edge-vars}). Since
	$d_{\mathbf p}=\sum_j|\tilde D_j|$ and $x\mapsto\sum_j|x_j|$ is continuous, the
	continuous mapping theorem gives \eqref{eq:weak-limit}. The normalized edge
	differences have uniformly bounded fourth moments, so
	$\{(d_{\mathbf p}/\sqrt N)^2\}_N$ is uniformly integrable and the second moment
	converges; evaluating
	$\Var(W)=\sum_j(1-\tfrac2\pi)\Var(Y_j)
	+\tfrac4\pi\sum_{j<k}\sqrt{\Var(Y_j)\Var(Y_k)}\,\Phi(\rho_{j,k})$ by the Nabeya
	identity reproduces \eqref{eq:Var-TVC} divided by $N$.
\end{proof}

	The proof of Theorem~\ref{thm:Var-TVC} is given after the following pair estimate.

\begin{proposition}[Bivariate $O(1)$ bound for an arbitrary pair of edges]
	\label{prop:pair-O1}
	Fix an interior probability vector $\mathbf p$ and two edges $j<k$ whose bins do not
	exhaust $\{1,\dots,n\}$ (automatic for $n\ge5$), put
	\[
		\mathcal B:=\{j,j+1,k,k+1\},
	\]
	and define
	\begin{equation}\label{eq:q0-def}
		q_0:=1-\sum_{i\in\mathcal B}p_i=\sum_{i\notin\mathcal B}p_i,
	\end{equation}
	where repeated indices in $\mathcal B$ are counted only once. Then $q_0>0$ and,
	as $N\to\infty$ at fixed $\mathbf p$,
	\begin{equation}\label{eq:pair-O1}
		\E\bigl[\,|\tilde D_j|\,|\tilde D_k|\,\bigr]
		=\frac2\pi\,\sigma_j\sigma_k
		\bigl[\sqrt{1-\rho_{j,k}^2}+\rho_{j,k}\arcsin\rho_{j,k}\bigr]+O(1),
	\end{equation}
	and consequently
	$\Cov(|\tilde D_j|,|\tilde D_k|)
	=\tfrac2\pi\,\sigma_j\sigma_k\,\Phi(\rho_{j,k})+O(1)$.
\end{proposition}

	Here $\mathcal B$ is the set of bins involved in the two edges: $|\mathcal B|=3$
	for an adjacent pair ($k=j+1$) and $|\mathcal B|=4$ for a disjoint pair
	($k\ge j+2$). The number $q_0$ is the total probability of the bins outside the
	pair. Its positivity is not an extra hypothesis in the present proposition:
	since $|\mathcal B|\le4<n$ and $\mathbf p$ is interior, at least one
	positive-probability bin remains outside $\mathcal B$. More generally, $q_0>0$
	for every pair that does not exhaust the bins; the only exceptions, the
	adjacent pair at $n=3$ and the disjoint pair at $n=4$, have $q_0=0$ and are
	treated separately in \S\ref{sec:small-n}.
	The constant supplied by the proof is uniform in $N$ but depends on the pair $(j,k)$
	and on $\mathbf p$ through the local-expansion constants of the increment law
	(aperiodicity, hence $q_0>0$, and nonsingularity); it is not uniform as
	$\min_i p_i\to0$ or $|\rho_{j,k}|\to1$, although the remainder itself shows no
	sensitivity to $q_0\to0$ or $|\rho_{j,k}|\to1$ in exact computations, and we make no
	uniform claim.

\medskip

\noindent\emph{Proof of Proposition~\ref{prop:pair-O1}.}\enspace
We prove the two displays through a sequence of claims. Abbreviate
$\rho:=\rho_{j,k}$,
$\delta_j:=\E D_j/N=p_{j+1}-p_j$, $\delta_k:=\E D_k/N=p_{k+1}-p_k$,
$c_j:=N\delta_j=\E D_j$, $c_k:=N\delta_k=\E D_k$, and use the standardized
coordinates $v=(v_1,v_2)$, $v_1=(y_1-c_j)/\sigma_j$, $v_2=(y_2-c_k)/\sigma_k$;
for a lattice point $m\in\Z^2$ write $v_m$ for its standardized position.
Represent the pair as a sum of i.i.d.\ increments,
\[
	(D_j,D_k)=\sum_{t=1}^N\zeta_t,
\]
where $\zeta_t$ records the change of $(D_j,D_k)$ produced by trial $t$. The
first four claims are properties of the single increment $\zeta$.

\begin{claim}[Increment law]\label{cl:incr}
	The distribution of $\zeta$ is, for a \emph{disjoint} pair ($|k-j|\ge2$),
	\[
		(-1,0),(+1,0),(0,-1),(0,+1),(0,0)\quad\text{with probabilities}\quad
		p_j,p_{j+1},p_k,p_{k+1},q_0,
	\]
	and, for an \emph{adjacent} pair ($k=j+1$, shared bin $j+1$),
	\[
		(-1,0),(+1,-1),(0,+1),(0,0)\quad\text{with probabilities}\quad
		p_j,p_{j+1},p_{j+2},q_0 .
	\]
	In both cases the probabilities sum to $1$.
\end{claim}
\begin{proof}
	A single trial falls in bin $i$ with probability $p_i$ and increments $X_i$ by
	$1$. In the disjoint case the four bins $j,j+1,k,k+1$ are distinct:
	incrementing $X_j$ changes $D_j=X_{j+1}-X_j$ by $-1$ and leaves $D_k$ fixed,
	giving $(-1,0)$; $X_{j+1}\to(+1,0)$; $X_k\to(0,-1)$; $X_{k+1}\to(0,+1)$; any of
	the remaining $n-4$ bins leaves both differences fixed, giving $(0,0)$ with the
	pooled probability $q_0=1-p_j-p_{j+1}-p_k-p_{k+1}$. In the adjacent case
	$k=j+1$ the bins are $j,j+1,j+2$: $X_j\to(-1,0)$; raising the \emph{shared}
	$X_{j+1}$ increases $D_j$ by $1$ and decreases $D_{j+1}=X_{j+2}-X_{j+1}$ by
	$1$, i.e.\ $(+1,-1)$; $X_{j+2}\to(0,+1)$; the other $n-3$ bins give $(0,0)$
	with $q_0=1-p_j-p_{j+1}-p_{j+2}$. In each case the listed events partition the
	$n$ bins, so the probabilities sum to $1$.
\end{proof}

\begin{claim}[Finiteness of moments]\label{cl:bdd}
	$\|\zeta\|_\infty\le1$; in particular $\zeta$ has finite moments of every order.
\end{claim}
\begin{proof}
	Every support point in Claim~\ref{cl:incr} has coordinates in $\{-1,0,1\}$, so
	$\zeta$ is bounded, and a bounded random vector has moments of all orders.
\end{proof}

\begin{claim}[Nonsingular per-trial covariance; strict correlation bound]
	\label{cl:nonsing}
	The per-trial covariance matrix $\Sigma:=\Cov(\zeta)$ is positive definite;
	consequently $|\rho_{j,k}|<1$ strictly, and $\rho_{j,k}$ does not depend on $N$.
\end{claim}
\begin{proof}
	In each case the support contains three affinely independent points of positive
	probability: $(0,0),(1,0),(0,1)$ in the disjoint case (probabilities
	$q_0,p_{j+1},p_{k+1}$, all positive), and $(0,0),(-1,0),(0,1)$ in the adjacent
	case (probabilities $q_0,p_j,p_{j+2}$, all positive). Suppose $u^\top\Sigma u=0$
	for a unit vector $u\in\R^2$. Then $\Var(u^\top\zeta)=0$, so $u^\top\zeta$ is
	almost surely equal to a constant $b$, i.e.\ every support point of $\zeta$ lies
	on the line $\{x:u^\top x=b\}$; but three affinely independent points cannot lie
	on one line. Hence $\Sigma\succ0$.

	Since $(D_j,D_k)=\sum_{t=1}^N\zeta_t$ is a sum of $N$ i.i.d.\ increments, its
	covariance matrix is $N\Sigma$; thus $\sigma_j^2=N\Sigma_{11}$,
	$\sigma_k^2=N\Sigma_{22}$, $\Cov(D_j,D_k)=N\Sigma_{12}$, and
	\[
		\rho_{j,k}=\frac{N\Sigma_{12}}{\sqrt{N\Sigma_{11}}\,\sqrt{N\Sigma_{22}}}
		=\frac{\Sigma_{12}}{\sqrt{\Sigma_{11}\Sigma_{22}}},
	\]
	independent of $N$. Positive definiteness gives
	$\Sigma_{11}\Sigma_{22}-\Sigma_{12}^2>0$, i.e.\ $|\rho_{j,k}|<1$ strictly.
	(Equivalently $\sigma_j^2/N=s_j-\delta_j^2\ge s_j(1-s_j)>0$ is a fixed positive
	constant, so $\sigma_j\asymp\sqrt N$.)
\end{proof}

\begin{claim}[Aperiodicity: the minimal lattice is $\Z^2$]\label{cl:aper}
	The support of $\zeta$ is contained in no proper coset $b+L'$ of a proper
	sublattice $L'\subsetneq\Z^2$; equivalently, the walk $(D_j,D_k)$ is strongly
	aperiodic on $\Z^2$.
\end{claim}
\begin{proof}
	Because $q_0>0$, the origin $(0,0)$ is a support point. Suppose the support lay
	in $b+L'$ for a sublattice $L'$ and some $b\in\Z^2$. From $(0,0)\in b+L'$ we get
	$-b\in L'$, hence $b\in L'$ and $b+L'=L'$, so the support lies in $L'$ itself. In
	the disjoint case $(1,0),(0,1)\in L'$; in the adjacent case $(-1,0),(0,1)\in L'$,
	and then $(1,0)=(1,-1)+(0,1)\in L'$ as well. In either case
	$L'\supseteq\Z(1,0)+\Z(0,1)=\Z^2$, forcing $L'=\Z^2$. Thus no proper coset
	contains the support, which is strong aperiodicity. (Had $q_0=0$ the origin need
	not lie in the support and a proper coset can occur: the disjoint pair that
	exhausts the four bins at $n=4$ satisfies $D_j+D_k\equiv N\pmod2$, a coset of the
	index-$2$ sublattice --- the reason such pairs are excluded here.)
\end{proof}

We now compare the lattice law with the two-dimensional Edgeworth approximant of
Bhattacharya--Rao~\cite[Theorem~22.1]{bhattacharya_rao},
written in the marginally standardized coordinates $v$
(equivalent to Mahalanobis coordinates since $|\rho|<1$). We
keep terms through order $N^{-2}$ --- the four correction terms of the $s=6$
instance used in Claim~\ref{cl:llt} ---
\begin{equation}\label{eq:fN}
	f_N(y)=\frac{\phi_\rho(v)}{\sigma_j\sigma_k}
	\Bigl[1+\frac{P_1(v)}{\sqrt N}+\frac{P_2(v)}{N}
	+\frac{P_3(v)}{N^{3/2}}+\frac{P_4(v)}{N^{2}}\Bigr],
\end{equation}
where $\phi_\rho$ is the standard bivariate normal density with correlation
$\rho$ and $P_r$ is the Cram\'er--Edgeworth polynomial of parity $(-1)^r$ --- so
$P_1,P_3$ are odd and $P_2,P_4$ even under $v\mapsto-v$ --- a fixed polynomial in
the \emph{per-trial} cumulants of $\zeta$, independent of $N$. Here $f_N$ is a
\emph{signed} approximant, not a probability density (the bracket may be
negative); it enters only integrated against the weight
$|y_1-c_j|\,|y_2-c_k|$ and through the lattice sum below.

\begin{claim}[Local limit bound]\label{cl:llt}
	Uniformly over $m\in\Z^2$,
	\begin{equation}\label{eq:llt}
		(1+|v_m|^6)\,\bigl|\Pr\{(D_j,D_k)=m\}-f_N(m)\bigr|=o(N^{-3}).
	\end{equation}
\end{claim}
\begin{proof}
	Apply the local Edgeworth expansion of
	Bhattacharya--Rao~\cite[Theorem~22.1]{bhattacharya_rao}
	with $d=2$ and $s=6$. Its hypotheses hold for
	$\zeta$: minimal lattice $\Z^2$ (Claim~\ref{cl:aper}), nonsingular covariance
	(Claim~\ref{cl:nonsing}), and finite absolute moments of every order since
	$\zeta$ is bounded (Claim~\ref{cl:bdd}), in particular of order $6$. The
	theorem, with $k=2$ and $s=6$, bounds the absolute lattice probabilities:
	$\sup_m(1+|v_m|^6)\,|\Pr\{(D_j,D_k)=m\}-f_N(m)|=o(N^{-(k+s-2)/2})=o(N^{-3})$ (its weight
	uses $(m-N\E\zeta)/\sqrt N$, equivalent to $v_m$ up to constants depending only on
	$\mathbf p$), after multiplying back by the cell factor $(\sigma_j\sigma_k)^{-1}=O(N^{-1})$.
	The $s-2=4$ correction terms of the approximant are exactly those of
	\eqref{eq:fN}; the first omitted term is of relative order
	$N^{-5/2}$. The approximant depends on $N$ only through the explicit scalings
	and the fixed polynomials $P_1,\dots,P_4$, so the implied constants are uniform
	in $N$.
\end{proof}

\begin{claim}[Passage to the continuous sum]\label{cl:riemann}
	With $h(y):=|y_1-c_j|\,|y_2-c_k|\,f_N(y)$,
	\begin{equation}\label{eq:step1}
		\E\bigl[|\tilde D_j|\,|\tilde D_k|\bigr]=\sum_{m\in\Z^2}h(m)+o(N^{-1}).
	\end{equation}
\end{claim}
\begin{proof}
	By definition
	$\E[|\tilde D_j||\tilde D_k|]
	=\sum_{m\in\Z^2}|m_1-c_j|\,|m_2-c_k|\,\Pr\{(D_j,D_k)=m\}$. Subtract
	$\sum_m h(m)$ and bound termwise. Writing
	$|m_1-c_j|\,|m_2-c_k|=\sigma_j\sigma_k|v_{m,1}|\,|v_{m,2}|$ and inserting the
	factor $(1+|v_m|^6)/(1+|v_m|^6)$,
	\begin{align*}
		&\sum_{m}\bigl|\,|m_1-c_j||m_2-c_k|\Pr\{(D_j,D_k)=m\}-h(m)\bigr|\\
		\qquad&=\sigma_j\sigma_k\sum_m g(v_m)\,(1+|v_m|^6)
		 \bigl|\Pr\{\cdot=m\}-f_N(m)\bigr|\\
		\qquad&\le \sigma_j\sigma_k\cdot o(N^{-3})\sum_m g(v_m),
	\end{align*}
	where $g(v):=|v_1v_2|/(1+|v|^6)$ and the last line uses \eqref{eq:llt}. Now $g$
	is nonnegative, continuous, and decays like $|v|^{-4}$ as $|v|\to\infty$, hence
	is integrable on $\R^2$ (this decay --- i.e.\ the weight exponent $6$ --- is
	exactly what makes $g$ integrable). The standardized lattice $\{v_m\}$ has
	spacings $\sigma_j^{-1}\times\sigma_k^{-1}$, i.e.\ cell area
	$(\sigma_j\sigma_k)^{-1}$; since $g$ is Riemann integrable with an integrable
	monotone tail,
	\[
		\sum_m g(v_m)=\sigma_j\sigma_k\!\iint_{\R^2}g\,dv\,(1+o(1))=O(\sigma_j\sigma_k).
	\]
	Therefore the total error is
	$\sigma_j\sigma_k\cdot o(N^{-3})\cdot O(\sigma_j\sigma_k)
	=o(N^{-3})\,O\bigl((\sigma_j\sigma_k)^2\bigr)=o(N^{-3})\,O(N^2)=o(N^{-1})$,
	using $\sigma_j\sigma_k=O(N)$ (Claim~\ref{cl:nonsing}). This is \eqref{eq:step1}.
\end{proof}

The sum $\sum_m h(m)$ is now evaluated by Poisson summation. Fix $N$; then
$h(y)=|y_1-c_j|\,|y_2-c_k|\,f_N(y)$ is continuous and decays faster than any
polynomial (Gaussian factor). Write
$\hat h(\xi):=\iint_{\R^2}h(y)\,e^{-i\xi\cdot y}\,dy$.

\begin{claim}[Tensor Fourier decay]\label{cl:decay}
	For each fixed $N$ there is $C_N<\infty$ with
	\begin{equation}\label{eq:decay}
		|\hat h(\xi)|\le C_N\,(1+\xi_1^2)^{-1}(1+\xi_2^2)^{-1},\qquad\xi\in\R^2 .
	\end{equation}
\end{claim}
\begin{proof}
	Consider the mixed distributional derivative
	$\partial_{y_1}^2\partial_{y_2}^2 h$. Applying the kink identity
	$\partial_y^2|y-c|=2\delta_c$ (Lemma~\ref{lem:imported}\textup{(ii)}) in each variable
	separately, this derivative is a \emph{finite signed measure}: one atom
	$4f_N(c_j,c_k)\,\delta_{(c_j,c_k)}$ (both kinks differentiated twice); two line
	measures on the kink lines $y_1=c_j$ and $y_2=c_k$ with rapidly decaying
	densities (one kink differentiated twice, the smooth factor in the other
	variable); and an absolutely continuous part
	$|y_1-c_j||y_2-c_k|\,\partial_{y_1}^2\partial_{y_2}^2 f_N$ together with the
	$\operatorname{sgn}$--$\partial$ cross terms, all integrable because $f_N$ is a
	Gaussian times a polynomial. Each piece has finite total variation, so
	$M_N:=\|\partial_{y_1}^2\partial_{y_2}^2 h\|_{TV}<\infty$. Taking Fourier
	transforms turns $\partial_{y_i}^2$ into multiplication by
	$(i\xi_i)^2=-\xi_i^2$, so $|\xi_1^2\xi_2^2\,\hat h(\xi)|\le M_N$. Combining this
	with the three cruder bounds obtained the same way (from $\|h\|_1$,
	$\|\partial_{y_1}^2h\|_{TV}$, $\|\partial_{y_2}^2h\|_{TV}$, all finite) gives
	$|\hat h(\xi)|\le C_N\min(1,\xi_1^{-2})\min(1,\xi_2^{-2})$, which is
	\eqref{eq:decay}. The decisive feature is that the bound \emph{factorizes} over
	the two frequencies.
\end{proof}

\begin{claim}[Validity of Poisson summation]\label{cl:psf}
	$\displaystyle\sum_{r\in\Z^2}|\hat h(2\pi r)|<\infty$, and
	\begin{equation}\label{eq:psf}
		\sum_{m\in\Z^2}h(m)=\sum_{r\in\Z^2}\hat h(2\pi r),
	\end{equation}
	with the frequency $2\pi r$ carried inside $\hat h$ (no factor $(2\pi)^{-2}$ on
	the right).
\end{claim}
\begin{proof}
	By \eqref{eq:decay},
	\[
		\sum_{r\in\Z^2}|\hat h(2\pi r)|
		\le C_N\Bigl(\sum_{r_1\in\Z}(1+4\pi^2r_1^2)^{-1}\Bigr)
		\Bigl(\sum_{r_2\in\Z}(1+4\pi^2r_2^2)^{-1}\Bigr)<\infty,
	\]
	each one-dimensional series converging and the tensor bound making the double
	series their product. (The isotropic bound $|\hat h(\xi)|\lesssim|\xi|^{-2}$
	would \emph{not} suffice, as $\sum_{r\neq0}|r|^{-2}$ diverges in $\Z^2$; the
	factorized form is what is needed.) Since $h$ is continuous, decays
	super-polynomially, and $\hat h$ is absolutely summable over the integer
	lattice, the standard Poisson summation theorem under absolute convergence
	\cite[Ch.~VII]{stein_weiss} applies: the periodization $\sum_m h(\cdot+m)$ is
	continuous with an absolutely convergent Fourier series and \eqref{eq:psf}
	holds; evaluating at $0$ gives the stated identity.
\end{proof}

\begin{claim}[Zero Fourier mode]\label{cl:zero}
	\begin{equation}\label{eq:zero}
		\hat h(0)=\frac2\pi\,\sigma_j\sigma_k\bigl[\sqrt{1-\rho^2}+\rho\arcsin\rho\bigr]
		+O(1).
	\end{equation}
\end{claim}
\begin{proof}
	By the substitution $y\mapsto v$ (Jacobian $\sigma_j\sigma_k$),
	\begin{equation}\label{eq:zero-int}
		\hat h(0)=\iint h\,dy
		=\sigma_j\sigma_k\!\iint_{\R^2}|v_1|\,|v_2|\,\phi_\rho(v)
		\Bigl[1+\frac{P_1(v)}{\sqrt N}+\frac{P_2(v)}{N}
		+\frac{P_3(v)}{N^{3/2}}+\frac{P_4(v)}{N^{2}}\Bigr]dv .
	\end{equation}
	Treat the leading term, the odd corrections and the even corrections in turn.

	\textit{Gaussian part.} $\sigma_j\sigma_k\iint|v_1||v_2|\phi_\rho(v)\,dv$. Here
	$\iint|v_1||v_2|\phi_\rho(v)\,dv=\E|Z_1||Z_2|$ for a standard bivariate normal
	$(Z_1,Z_2)$ of correlation $\rho$, which by the Nabeya identity
	\eqref{eq:sheppard-mardia} equals
	$\frac2\pi[\sqrt{1-\rho^2}+\rho\arcsin\rho]$; the
	part is thus $\frac2\pi\sigma_j\sigma_k[\sqrt{1-\rho^2}+\rho\arcsin\rho]$.

	\textit{The odd corrections vanish identically.} For $\ell\in\{1,3\}$ let
	\[
		I_\ell:=\iint_{\R^2}|v_1||v_2|\,\phi_\rho(v)\,P_\ell(v)\,dv.
	\]
	Under $v\mapsto-v$ the weight
	$|v_1||v_2|$ is even and $\phi_\rho(-v)=\phi_\rho(v)$ (the density depends on $v$
	only through the even quadratic form $v^\top\Sigma_\rho^{-1}v$), while $P_\ell$
	is odd; the integrand is odd and $I_\ell=-I_\ell=0$. Thus the $N^{-1/2}$ and
	$N^{-3/2}$ parts both vanish. This parity cancellation is why no $O(\sqrt N)$
	term appears in \eqref{eq:pair-O1}, although the third cumulants of a
	non-uniform multinomial do not vanish.

	\textit{The even corrections are $O(1)$.} For $\ell\in\{2,4\}$ the coefficients
	$I_\ell:=\iint|v_1||v_2|\phi_\rho(v)P_\ell(v)\,dv$ are finite constants
	($P_\ell$ a fixed polynomial; all Gaussian moments of
	$|v_1||v_2|\cdot(\text{polynomial})$ are finite). Hence the $N^{-1}$ part is
	$\frac{\sigma_j\sigma_k}{N}I_2=O(\sigma_j\sigma_k/N)=O(1)$ and the $N^{-2}$ part
	is $\frac{\sigma_j\sigma_k}{N^{2}}I_4=O(\sigma_j\sigma_k/N^{2})=O(N^{-1})=o(1)$.
	(For a disjoint pair the $\kappa_{2,2}$ portion of $I_2$ is exactly the explicit
	term of Proposition~\ref{prop:disj-dominant}; we need only that it is $O(1)$.)
	Adding the parts gives \eqref{eq:zero}.
\end{proof}

For the nonzero modes we separate \emph{axis} modes ($r_1r_2=0$) from
\emph{double} modes ($r_1r_2\neq0$).

\begin{claim}[Axis modes]\label{cl:axis}
	For $r=(r_1,0)$, $r_1\neq0$,
	\begin{equation}\label{eq:axis}
		\hat h(2\pi r_1,0)=-\frac{2\,g(c_j)\,e^{-2\pi ir_1c_j}}{(2\pi r_1)^2}
		+O\!\bigl(N^{-1}|r_1|^{-3}\bigr),\qquad g(c_j)=O(1),
	\end{equation}
	and $\sum_{r_1\neq0}\hat h(2\pi r_1,0)=O(1)$; the family $(0,r_2)$ is symmetric.
\end{claim}
\begin{proof}
	Integrating out $y_2$ first,
	\[
		\hat h(2\pi r_1,0)=\int_\R|y_1-c_j|\,g(y_1)\,e^{-2\pi ir_1y_1}\,dy_1,
		\qquad g(y_1):=\int_\R|y_2-c_k|\,f_N(y_1,y_2)\,dy_2 .
	\]
	With $t=(y_1-c_j)/\sigma_j$ and the substitution $y_2\mapsto v_2$,
	\[
		g(y_1)=\frac{\sigma_k}{\sigma_j}\int_\R|v_2|\,\phi_\rho(t,v_2)
		\Bigl[\sum_{\ell=0}^{4}N^{-\ell/2}P_\ell(t,v_2)\Bigr]dv_2
		=\frac{\sigma_k}{\sigma_j}\sum_{\ell=0}^{4}N^{-\ell/2}H_\ell(t),
	\]
	where $P_0\equiv1$ and $H_\ell(t):=\int|v_2|\phi_\rho(t,v_2)P_\ell(t,v_2)\,dv_2$.
	Each $H_\ell$ is a Schwartz function of $t$ (a Gaussian in $t$ times a
	polynomial), the family bounded uniformly in $N$ since $P_0,\dots,P_4$ are
	$N$-independent; moreover $H_0,H_2,H_4$ are even and $H_1,H_3$ odd. Indeed,
	under $t\mapsto-t$ substitute $v_2\mapsto-v_2$: $|v_2|$ is unchanged,
	$\phi_\rho(-t,-v_2)=\phi_\rho(t,v_2)$, and $P_\ell(-t,-v_2)=(-1)^\ell
	P_\ell(t,v_2)$, giving $H_\ell(-t)=(-1)^\ell H_\ell(t)$.

	Thus $g(y_1)=A_N\,H\bigl((y_1-c_j)/\sigma_j\bigr)$ with $A_N=\sigma_k/\sigma_j$,
	$\sigma=\sigma_j\asymp\sqrt N$, and $H=\sum_{\ell=0}^{4}N^{-\ell/2}H_\ell$ a
	bounded Schwartz family. The scaled-profile Fourier estimate
	(Lemma~\ref{lem:imported}\textup{(iii)}), with $c=c_j$ and $\omega=2\pi r_1$, gives
	\[
		\hat h(2\pi r_1,0)
		=-\frac{2\,g(c_j)\,e^{-2\pi ir_1c_j}}{(2\pi r_1)^2}+O\!\bigl(N^{-1}|r_1|^{-3}\bigr).
	\]
	At the kink $t=0$ the odd parts vanish ($H_1(0)=H_3(0)=0$), so
	$g(c_j)=\frac{\sigma_k}{\sigma_j}[H_0(0)+N^{-1}H_2(0)+N^{-2}H_4(0)]
	=O(\sigma_k/\sigma_j)=O(1)$,
	as $\sigma_k/\sigma_j$ is a fixed positive constant. Summing, the leading terms
	are $O(r_1^{-2})$ and the remainders $O(N^{-1}|r_1|^{-3})$; both series
	converge, so $\sum_{r_1\neq0}\hat h(2\pi r_1,0)=O(1)$. The modes $(0,r_2)$ are
	handled identically by integrating out $y_1$ first. The phases
	$e^{-2\pi ir_1c_j}=e^{-2\pi ir_1N\delta_j}$ are the two-dimensional lattice
	oscillation (the analogue of $B_2(\{N\delta_j\})$) and are absorbed into the $O(1)$
	remainder.
\end{proof}

\begin{claim}[Double modes]\label{cl:double}
	$\displaystyle\sum_{r_1r_2\neq0}\hat h(2\pi r)=O(N^{-1})$.
\end{claim}
\begin{proof}
	For $r_1r_2\neq0$ apply the kink identity in \emph{both} variables. The tensor
	product of the two leading kink terms gives
	\[
		\hat h(2\pi r)=\frac{4\,f_N(c_j,c_k)\,e^{-2\pi i(r_1c_j+r_2c_k)}}
		{(2\pi r_1)^2(2\pi r_2)^2}+R_N(r),
	\]
	where
	\[
		|R_N(r)|\le C N^{-1}\bigl(|r_1|^{-3}|r_2|^{-2}
		+|r_1|^{-2}|r_2|^{-3}+|r_1|^{-3}|r_2|^{-3}\bigr).
	\]
	Indeed, after one kink is extracted in one variable, the remaining
	one-dimensional profile in the other variable has the form covered by
	Lemma~\ref{lem:imported}\textup{(iii)}. The even Edgeworth pieces have zero
	linear term at the kink, while the odd pieces carry the explicit factors
	$N^{-1/2}$ and $N^{-3/2}$, so the derivative condition
	$H_N'(0)=O(N^{-1/2})$ holds. Pairing one leading kink term with one such
	remainder gives the first two displayed bounds for $R_N$; pairing two
	remainders gives the third. By \eqref{eq:fN},
	\[
		f_N(c_j,c_k)=\frac{\phi_\rho(0)}{\sigma_j\sigma_k}\bigl[1+O(N^{-1})\bigr]
		=O\bigl((\sigma_j\sigma_k)^{-1}\bigr)=O(N^{-1}),
	\]
	with $\phi_\rho(0)=(2\pi\sqrt{1-\rho^2})^{-1}$; the odd terms vanish at
	$v=0$. Hence $|\hat h(2\pi r)|$ is bounded by a constant times the summable
	envelope
	\[
		N^{-1}\bigl(|r_1r_2|^{-2}+|r_1|^{-3}|r_2|^{-2}
		+|r_1|^{-2}|r_2|^{-3}+|r_1r_2|^{-3}\bigr),
	\]
	and the whole double-mode family sums to $O(N^{-1})$.
\end{proof}

	\paragraph{A direct bound} The axis- and double-mode orders can also be read off without
	the imported estimate~\textup{(iii)}: writing the joint folded density in standardized
	coordinates, $h(y)=K_N(v)$ with
	$K_N(v)=|v_1||v_2|\phi_\rho(v)\bigl[1+\sum_{\ell\ge1}N^{-\ell/2}P_\ell(v)\bigr]$ a family
	bounded in $N$ together with its mixed distributional derivatives, one has
	$\|\partial_{y_1}^2\partial_{y_2}^2h\|_{\mathrm{TV}}=(\sigma_j\sigma_k)^{-1}
	\|\partial_{v_1}^2\partial_{v_2}^2K_N\|_{\mathrm{TV}}=O(N^{-1})$ and
	$\|\partial_{y_1}^2h\|_{\mathrm{TV}}=O(1)$, so the double modes contribute $O(N^{-1})$ and the
	axis modes $O(1)$ uniformly in $N$; the explicit axis-mode computation below identifies the
	oscillation.

\begin{claim}[All nonzero modes]\label{cl:nonzero}
	$\displaystyle\sum_{r\neq0}\hat h(2\pi r)=O(1)$.
\end{claim}
\begin{proof}
	By Claims~\ref{cl:axis} and~\ref{cl:double} the axis families contribute $O(1)$
	and the double family $O(N^{-1})=o(1)$; the total over $r\neq0$ is $O(1)$, the
	rearrangement being justified by the absolute convergence of
	Claim~\ref{cl:psf}.
\end{proof}

\begin{claim}[Product form]\label{cl:main}
	$\displaystyle\E\bigl[|\tilde D_j|\,|\tilde D_k|\bigr]
	=\frac2\pi\sigma_j\sigma_k\bigl[\sqrt{1-\rho^2}+\rho\arcsin\rho\bigr]+O(1)$.
\end{claim}
\begin{proof}
	By Poisson summation \eqref{eq:psf},
	$\sum_m h(m)=\hat h(0)+\sum_{r\neq0}\hat h(2\pi r)$; inserting the zero mode
	\eqref{eq:zero} and Claim~\ref{cl:nonzero},
	$\sum_m h(m)=\frac2\pi\sigma_j\sigma_k[\sqrt{1-\rho^2}+\rho\arcsin\rho]+O(1)$.
	Combining with \eqref{eq:step1} (Claim~\ref{cl:riemann}) gives the assertion,
	the $o(1)$ being absorbed in the $O(1)$.
\end{proof}

\begin{claim}[Covariance form]\label{cl:cov}
	$\displaystyle\Cov\bigl(|\tilde D_j|,|\tilde D_k|\bigr)
	=\frac2\pi\sigma_j\sigma_k\,\Phi(\rho)+O(1)$.
\end{claim}
\begin{proof}
	By the per-edge folded mean (Lemma~\ref{lem:imported}\textup{(i)}),
	$\E|\tilde D_j|=\sqrt{2/\pi}\,\sigma_j(1+O(N^{-1}))$ and likewise for $k$, so
	\[
		\E|\tilde D_j|\,\E|\tilde D_k|
		=\frac2\pi\sigma_j\sigma_k\bigl(1+O(N^{-1})\bigr)=\frac2\pi\sigma_j\sigma_k+O(1),
	\]
	the error being $\frac2\pi\sigma_j\sigma_k\cdot O(N^{-1})=O(1)$. Subtracting from
	Claim~\ref{cl:main},
	\[
		\Cov\bigl(|\tilde D_j|,|\tilde D_k|\bigr)
		=\frac2\pi\sigma_j\sigma_k\bigl[\sqrt{1-\rho^2}+\rho\arcsin\rho-1\bigr]+O(1)
		=\frac2\pi\sigma_j\sigma_k\,\Phi(\rho)+O(1),
	\]
	by the definition $\Phi(\rho)=\sqrt{1-\rho^2}+\rho\arcsin\rho-1$ of
	\eqref{eq:SM-Phi}.
\end{proof}

\noindent
Claims~\ref{cl:main} and~\ref{cl:cov} are the two displays of the proposition,
so the proof of Proposition~\ref{prop:pair-O1} is complete.\hfill$\qed$

\medskip

We can now assemble the theorem.

\begin{proof}[Proof of Theorem~\ref{thm:Var-TVC}]
	By definition,
	\[
		\Var(d_{\mathbf p})=\sum_{j=1}^{n-1}\Var(|\tilde D_j|)
		+2\sum_{1\le j<k\le n-1}\Cov(|\tilde D_j|,|\tilde D_k|).
	\]

	\emph{Per-edge term.} Exactly $\E\tilde D_j^2=\Var(D_j)=\sigma_j^2$, while
	$\E|\tilde D_j|=\sqrt{2/\pi}\,\sigma_j(1+O(N^{-1}))$ by the per-edge expansion
	(Lemma~\ref{lem:imported}\textup{(i)}); hence
	$\Var(|\tilde D_j|)=\sigma_j^2-\tfrac2\pi\sigma_j^2(1+O(N^{-1}))
	=(1-2/\pi)\sigma_j^2+O(1)$, the Gaussian value $\Var|Z|=1-2/\pi$ appearing as
	the leading coefficient.

	\emph{Covariance term.} The pair $(\tilde D_j/\sigma_j,\tilde D_k/\sigma_k)$ is
	asymptotically bivariate normal with correlation $\rho_{j,k}$, whose absolute
	covariance is $\frac2\pi\sigma_j\sigma_k\Phi(\rho_{j,k})$
	by the Nabeya identity \eqref{eq:sheppard-mardia}.
	Proposition~\ref{prop:pair-O1} shows that for
	$n\ge5$ this value is attained with an absolute $O(1)$ remainder,
	\[
		\Cov(|\tilde D_j|,|\tilde D_k|)=\frac2\pi\sigma_j\sigma_k\Phi(\rho_{j,k})+O(1).
	\]
	This holds for every pair $j<k$, adjacent and disjoint alike; only the
	magnitude of $\rho_{j,k}$ differs (for adjacent pairs $\rho_j=\rho_{j,j+1}$ is
	$-1/2$ in the uniform case and of either sign in general; for disjoint pairs
	$\rho_{j,k}=-(p_{j+1}-p_j)(p_{k+1}-p_k)/\sqrt{[s_j-\delta_j^2][s_k-\delta_k^2]}$,
	with $s_j=p_j+p_{j+1}$, $\delta_j=p_{j+1}-p_j$).

	\emph{Summation.} Combining the two contributions yields \eqref{eq:Var-TVC};
	the remainder $R_N$ collects the sub-leading $O(1)$ per edge and the finite-$N$
	Edgeworth correction, which is $O(1)$ per pair
	(Proposition~\ref{prop:pair-O1}).
\end{proof}

\begin{remark}[The remainder, the hypothesis $n\ge5$, and the regimes]
	\label{rem:VarTVC-scope}
	The sum over pairs in \eqref{eq:Var-TVC} splits into the \emph{adjacent} pairs
	$k=j+1$, with $\rho_{j,j+1}=\rho_j$ from \eqref{eq:rho-j-mult} (the closed form
	being \eqref{eq:Cov-Dj-adj} divided by $\sigma_j\sigma_{j+1}$), and the
	\emph{disjoint} pairs $k\ge j+2$, with $\rho_{j,k}$ from \eqref{eq:Cov-Dj-disj}
	divided by $\sigma_j\sigma_k$. The per-edge part of $R_N$ comes from
	\cite[Lemma~4.1]{banic_elezovic_p3}; every covariance pair --- adjacent and
	disjoint alike --- is covered by Proposition~\ref{prop:pair-O1}, whose hypothesis
	that no pair of edges exhausts the bins is guaranteed by $n\ge5$
	(Remark~\ref{rem:pair-boundary} and \S\ref{sec:small-n} treat $n=3,4$). The constant in
	$R_N=O(1)$ depends on $(n,\mathbf p)$ and is not uniform as $\min_jp_j\to0$ or
	as $\max_{j<k}|\rho_{j,k}|\to1$ (cf.\ the comb near $\varepsilon=1/n$ in
	\S\ref{ssec:comb-closed-form}). A growing-$n$ version would require uniform
	estimates not established here; accordingly, all large-$n$ statements later in
	the paper (Remarks~\ref{rem:disj-scaling} and~\ref{rem:comb-leading},
	Proposition~\ref{prop:smooth-universal}) concern the \emph{iterated} limit
	$N\to\infty$ followed by $n\to\infty$ of the leading coefficient
	$v(\mathbf p)$, not a joint regime.
\end{remark}

\begin{remark}[Boundary cases $n=3,4$]\label{rem:pair-boundary}
	The hypothesis $q_0>0$ excludes exactly two configurations: the adjacent pair
	at $n=3$ (the three bins exhaust, and the increment lattice degenerates to an
	index-$3$ sublattice) and the disjoint pair
	at $n=4$ (the four bins exhaust, $D_j+D_k\equiv N\pmod2$). In the uniform
	case the former is covered by the conditional analysis of the adjacent skeleton
	$\hadj$ in \cite[Lemma~5.1 and Theorem~5.3]{banic_elezovic_p2}, while the
	latter is covered by the exact $h(N)$
	reduction in the small-$n$ section of \cite{banic_elezovic_p2}
	(cf.\ \cite[Theorem~4.1]{banic_elezovic_p2} for the skeleton $h$). For general
	$\mathbf p$,
	both boundary configurations fall outside Proposition~\ref{prop:pair-O1}; they
	are settled in \S\ref{sec:small-n} (Theorem~\ref{thm:sublattice}) by Poisson
	summation over the sublattice coset, with the conclusion of
	Proposition~\ref{prop:pair-O1} intact. For $n\ge5$ no pair is affected.
\end{remark}

\noindent
We next isolate the part of the $O(1)$ remainder which survives in the uniform
disjoint covariance. This refinement is not needed for the leading variance
formula, but it explains how the present expansion connects with the uniform
theory of Part~2 and identifies one explicit Edgeworth contribution inside the
otherwise implicit remainder.

\begin{lemma}[Exact reduction of the disjoint covariance]\label{lem:disj-reduction}
	For a disjoint pair $|k-j|\ge2$ with $s_j:=p_j+p_{j+1}>0$ and
	$s_k:=p_k+p_{k+1}>0$, write the block sums $M_j := X_j + X_{j+1}$,
	$M_k := X_k + X_{k+1}$ and the conditional folded means
	\[
		\hat G_j(a) := \E\bigl[\,|\tilde D_j|\ \big|\ M_j = a\bigr]
		= \E\bigl|\,2\Bin(a,\pi_j) - a - N\delta_j\bigr|,
	\]
	\[
		\pi_j := \frac{p_{j+1}}{p_j+p_{j+1}},\ \ \delta_j := p_{j+1}-p_j,
	\]
	and $\hat G_k$ analogously. Then
	\begin{gather*}
		\Cov\bigl(|\tilde D_j|,|\tilde D_k|\bigr)
		=\Cov\bigl(\hat G_j(M_j),\hat G_k(M_k)\bigr),\\
		(M_j,M_k)\sim\Mult\bigl(N;\,s_j,\,s_k,\,1-s_j-s_k\bigr).
	\end{gather*}
\end{lemma}

\begin{proof}
	Condition on the block totals $(M_j,M_k)$. The two-block conditional
	independence is checked directly from the four-bin marginal
	(Proposition~\ref{prop:k-marginal}): for $x+y=a$, $u+v=b$, $a+b\le N$,
	\begin{multline*}
		\Pr\bigl(X_j=x,X_{j+1}=y,X_k=u,X_{k+1}=v\ \big|\ M_j=a,\,M_k=b\bigr)\\
		=\binom{a}{y}\pi_j^{\,y}(1-\pi_j)^{x}\;\binom{b}{v}\pi_k^{\,v}(1-\pi_k)^{u},
	\end{multline*}
	because in the quotient of the two multinomial pmfs the factors
	$(1-s_j-s_k)^{N-a-b}$ and the multinomial normalizations cancel and the
	remaining expression factors over the two blocks --- this is the two-block
	extension of Proposition~\ref{prop:2bin-factor}. Hence, given the totals, the
	intra-block counts
	$X_{j+1}\sim\Bin(M_j,\pi_j)$ and $X_{k+1}\sim\Bin(M_k,\pi_k)$ are independent
	and each
	depends only on its own total; so $\tilde D_j \perp \tilde D_k \mid (M_j,M_k)$ and
	$\E[\,|\tilde D_j| \mid M_j,M_k] = \hat G_j(M_j)$. Therefore
	$\E[\,|\tilde D_j|\,|\tilde D_k|\,] = \E[\hat G_j(M_j)\hat G_k(M_k)]$ and
	$\E|\tilde D_j| = \E[\hat G_j(M_j)]$; subtracting gives the identity. The pair
	$(M_j,M_k)$ is a two-coordinate marginal of the multinomial, hence trinomial with the
	stated parameters.
\end{proof}

\begin{remark}[Disjoint per-pair remainder]\label{prop:disj-O1}
	Let $\mathbf p$ be interior, let $n\ge5$, and let $|k-j|\ge2$. Then, as
	$N\to\infty$,
	\[
		\Cov\bigl(|\tilde D_j|,|\tilde D_k|\bigr)
		= \frac{2}{\pi}\sigma_j\sigma_k\,\Phi(\rho_{j,k}) + O(1),
	\]
	with the constant depending on $(n,\mathbf p,j,k)$, so each disjoint pair
	contributes $O(1)$ to the remainder $R_N$ in \eqref{eq:Var-TVC}. This is the
	special case of Proposition~\ref{prop:pair-O1} in which the two edges are
	disjoint: the pair uses four bins and, since $n\ge5$ with interior $\mathbf p$,
	the complement $q_0=\sum_{i\notin\{j,j+1,k,k+1\}}p_i>0$ required there is
	automatic. It is recorded separately because disjoint pairs are numerous in the
	variance sum, and because Lemma~\ref{lem:disj-reduction} gives the finite-$N$
	reduction from which the explicit $\kappa_{2,2}$ contribution of
	Proposition~\ref{prop:disj-dominant} is extracted.
\end{remark}

\begin{proposition}[The $\kappa_{2,2}$ contribution to the disjoint $O(1)$ remainder]\label{prop:disj-dominant}
	Let $\mathbf p$ be interior, let $|k-j|\ge2$, and assume
	$q_0:=1-s_j-s_k>0$. Put
	\[
		s_j:=p_j+p_{j+1},\quad s_k:=p_k+p_{k+1},\quad
		\delta_j:=p_{j+1}-p_j,\quad \delta_k:=p_{k+1}-p_k,
	\]
	\[
		\tilde s_j^2:=\Var(D_j)/N=s_j-\delta_j^2,\quad
		\tilde s_k^2:=\Var(D_k)/N=s_k-\delta_k^2,
		\quad
		\rho_{j,k}:=-\frac{\delta_j\delta_k}{\tilde s_j\tilde s_k}.
	\]
	Define the $O(1)$ remainder of Remark~\ref{prop:disj-O1} by
	\[
		C_{j,k}(N) := \Cov\bigl(|\tilde D_j|,|\tilde D_k|\bigr) - \tfrac2\pi\sigma_j\sigma_k\Phi(\rho_{j,k}),
	\]
	and let
	\[
		\kappa_{2,2}:=-s_js_k+2s_k\delta_j^2+2s_j\delta_k^2-6\delta_j^2\delta_k^2.
	\]
	Then
	\begin{equation}\label{eq:disj-dominant}
		C_{j,k}(N)
		= \frac{\kappa_{2,2}}{2\pi\,\tilde s_j\tilde s_k\sqrt{1-\rho_{j,k}^2}}
		+ C_{j,k}'(N),
	\end{equation}
	where the first term is the contribution of the per-trial fourth cumulant $\kappa_{2,2}$
	in the bivariate Edgeworth expansion of the joint law of $(\tilde D_j,\tilde D_k)$, and
	$C_{j,k}'(N)$ collects the remaining Edgeworth terms. At the exact uniform point
	$\delta_j=\delta_k=0$, $\rho_{j,k}=0$, $n\ge5$, one has
	\[
		\frac{\kappa_{2,2}}{2\pi\,\tilde s_j\tilde s_k\sqrt{1-\rho_{j,k}^2}}
		=-\frac1{\pi n},
		\qquad C_{j,k}'(N)=O(N^{-1}).
	\]
\end{proposition}

\noindent
Here $\kappa_{2,2}$ is the per-trial joint fourth cumulant of the two disjoint
edge increments. This generalizes the uniform computation
\cite[Proposition~4.5 and \S4]{banic_elezovic_p2}; note that $\kappa_{2,2}$ there is the
total cumulant (carrying the factor $N$), here the per-trial one. The term displayed in
\eqref{eq:disj-dominant} is only the explicit $\kappa_{2,2}$ part of the $O(1)$ correction. The remainder
$C_{j,k}'(N)$ collects the other joint fourth-cumulant portions
($\kappa_{3,1}$, $\kappa_{1,3}$), the marginal fourth-cumulant and
squared-third-cumulant terms of the bivariate Edgeworth polynomial, and the
two-dimensional lattice corrections. Thus the proposition does not claim that
the displayed term is dominant away from uniformity. At the uniform point the
extra portions cancel, and the last display matches, both in value and in rate,
the uniform disjoint theorem \cite[Theorem~4.6]{banic_elezovic_p2}.

\begin{proof}
	Per trial the edge increments $(Y_j,Y_k)$ satisfy $Y_jY_k\equiv0$ (disjoint edges: at most one of
	the bins $j,j{+}1,k,k{+}1$ is hit), so all mixed raw moments vanish and the moment generating
	function is additive,
	\[
		\E\,e^{uY_j+vY_k} = A(u)+B(v)+p_0,\qquad A(u) = p_je^{-u}+p_{j+1}e^{u},\ \ B(v) = p_ke^{-v}+p_{k+1}e^{v},
	\]
	$p_0 = 1-s_j-s_k$. With $S = A+B+p_0$, $K = \log S$, and
	$A(0)=A''(0)=s_j$, $A'(0)=\delta_j$
	(and likewise for $B$, while $S(0,0)=1$), first differentiate once in
	each variable:
	\[
		K_{uv}=-\frac{S_uS_v}{S^2}.
	\]
	Because $S_u$ depends only on $u$ and $S_v$ only on $v$, two further
	differentiations give
	\[
		K_{uuvv}
		=-\frac{S_{uu}S_{vv}}{S^2}
		+\frac{2S_{uu}S_v^2+2S_u^2S_{vv}}{S^3}
		-\frac{6S_u^2S_v^2}{S^4}.
	\]
	Evaluation at the origin therefore yields
	\[
		K_{uv}(0)=-\delta_j\delta_k,\qquad
		K_{uuvv}(0)
		=-s_js_k+2s_k\delta_j^2+2s_j\delta_k^2-6\delta_j^2\delta_k^2
		=\kappa_{2,2},
	\]
	the joint second and fourth per-trial cumulants. By additivity the joint fourth cumulant of
	$(\tilde D_j,\tilde D_k)$ is $N\kappa_{2,2}$, while $\rho_{j,k}$ is $N$-independent. Let
	$\phi_\Sigma$ be the centered bivariate Gaussian with covariance
	$\Sigma = \bigl(\begin{smallmatrix}\sigma_j^2 & \Cov(D_j,D_k)\\ \Cov(D_j,D_k) & \sigma_k^2\end{smallmatrix}\bigr)$,
	$\det\Sigma = \sigma_j^2\sigma_k^2(1-\rho_{j,k}^2)$. In the bivariate Edgeworth expansion of the
	joint density the term carrying $\kappa_{2,2}$ is
	$\tfrac14 N\kappa_{2,2}\,\partial_x^2\partial_y^2\phi_\Sigma$. Being a \emph{joint} cumulant it does
	not enter $\E|\tilde D_j|\,\E|\tilde D_k|$, so its whole contribution falls on the covariance:
	\[
		C_{j,k}^{(\kappa_{2,2})} = \frac{N\kappa_{2,2}}{4}\int_{\mathbb R^2}|x|\,|y|\,
		\partial_x^2\partial_y^2\phi_\Sigma(x,y)\,dx\,dy = N\kappa_{2,2}\,\phi_\Sigma(0,0),
	\]
	the kink integral $\int|x||y|\partial_x^2\partial_y^2\phi_\Sigma = 4\phi_\Sigma(0,0)$ being evaluated
	by two integrations by parts in each variable via $\int_{\mathbb R}|t|\,u''(t)\,dt = 2u(0)$. With
	$\phi_\Sigma(0,0) = 1/(2\pi\sqrt{\det\Sigma}) = 1/(2\pi\sigma_j\sigma_k\sqrt{1-\rho_{j,k}^2})$ and
	$\sigma_j\sigma_k = N\tilde s_j\tilde s_k$ this is the displayed
	joint-fourth-cumulant term in (\ref{eq:disj-dominant}). The
	integration of the non-smooth weight $|x||y|$ against the Edgeworth density is
	justified by the
	regularization argument of \cite[Remark~3.2]{banic_elezovic_p2}. The remaining
	Edgeworth portions listed in the statement, together with the two-dimensional
	lattice oscillation in $\{N\delta_j\},\{N\delta_k\}$, make up $C_{j,k}'(N)$; the
	uniform reduction is
	obtained as follows. Setting $\delta_j=\delta_k=0$ gives
	$\rho_{j,k}=0$, $\kappa_{2,2}=-s_js_k=-4/n^2$, and
	$\tilde s_j\tilde s_k=2/n$, so the displayed contribution is
	$-1/(\pi n)$; moreover at this point the joint odd cumulants vanish
	($\kappa_{3,1}=\kappa_{1,3}=0$) and the marginal Edgeworth terms factor and
	cancel exactly. The sharper residual estimate
	$C_{j,k}'(N)=O(N^{-1})$ is then exactly the content of
	\cite[Theorem~4.6]{banic_elezovic_p2}.
\end{proof}

\begin{remark}[Numerical confirmation]\label{rem:disj-dominant-num}
	The exact trinomial reduction (Lemma~\ref{lem:disj-reduction}) confirms
	(\ref{eq:disj-dominant}): for the bin masses
	$(p_1,p_2,p_3,p_4)=(0.1,0.3,0.25,0.15)$,
	with the remaining mass $0.2$ outside the pair, the exactly computed remainder
	for the disjoint pair $(1,3)$ is $C_{1,3}(N)=-0.0566$ at $N=100$ and
	$-0.0563$ at $N=300$, against the
	$\kappa_{2,2}$ term $-0.0521$; the gap
	$C_{1,3}'(N)\approx -0.004$ is nearly constant in $N$ for these masses
	(the oscillatory lattice part, of periods $5$ and $10$ in $N$ through
	$\{N\delta_j\},\{N\delta_k\}$, has peak-to-peak amplitude about $2\cdot10^{-4}$ near
	$N=100$ and $10^{-4}$ near $N=300$, decaying with $N$). At the exact
	uniform point the residual tends to zero and
	the displayed term is exactly $-1/(\pi n)$. The example suggests that the
	residual can be small near uniformity, but no uniform small-nonuniformity claim
	is made here, because for other masses the lattice part may oscillate more
	strongly. A closed form for this residual is not pursued.
\end{remark}

\begin{remark}[Reduction to the uniform disjoint theorem]\label{rem:disj-reduction-uniform}
	Under the uniform $p_i\equiv1/n$ one has $\delta_j=0$ and $\pi_j=\tfrac12$, so
	$\hat G_j(a) = \E|2\Bin(a,\tfrac12)-a| = g(a)$ is the symmetric-walk mean of
	the disjoint-skeleton analysis \cite[Theorem~4.1]{banic_elezovic_p2},
	$\rho_{j,k}=0$, and Lemma~\ref{lem:disj-reduction} becomes
	$\Cov(|D_j|,|D_k|) = \Cov(g(M_j),g(M_k))$ --- the representation underlying
	the uniform disjoint theorem \cite[Theorem~4.6]{banic_elezovic_p2}, whose limit
	is $-1/(\pi n)$ (with Gaussian part zero).
\end{remark}

\begin{remark}[The leading disjoint contribution vanishes for the uniform case]\label{rem:disj-uniform}
	For the uniform distribution, $p_{i+1} - p_i \equiv 0$, so $\rho_{j,k} = 0$ for all
	$|k-j| \ge 2$, and $\Phi(0) = 0$. Thus the \emph{leading Gaussian disjoint
	terms} in (\ref{eq:Var-TVC}) vanish. The exact disjoint covariances do not:
	they contribute the $O(1)$ residual $-1/(\pi n)+O(1/N)$ per pair, which is
	absorbed into $R_N$; see \cite[Theorem~4.6]{banic_elezovic_p2}.
\end{remark}

\begin{remark}[Leading scaling of the disjoint terms]\label{rem:disj-scaling}
	The scaling depends jointly on the local masses
	$s_j=p_j+p_{j+1}$ and gradients $\delta_j=p_{j+1}-p_j$:
	\[
		\rho_{j,k}=-\frac{\delta_j\delta_k}
		{\sqrt{(s_j-\delta_j^2)(s_k-\delta_k^2)}}.
	\]
	When the correlations are small, $\Phi(\rho)=\rho^2/2+O(\rho^4)$, so the
	disjoint contribution can be estimated directly from these two quantities.
	There is no universal growing-$n$ order without specifying a normalized family
	of probability vectors. In particular, for the symmetric comb
	$p_j=(1+(-1)^j a)/n$ with fixed $|a|<1$, one has
	$|\rho_{j,k}|=O(1/n)$, hence the total disjoint contribution is $O(N/n)$,
	not a leading $O(nN)$ term.
\end{remark}

	\paragraph{Properties of the function $\Phi$}
	The Nabeya correction
	$\Phi(\rho)=\sqrt{1-\rho^2}+\rho\arcsin\rho-1$ is even, strictly convex on
	$(-1,1)$, and nonnegative, with equality only at $\rho=0$. Moreover
	$\Phi(\rho)=\rho^2/2+O(\rho^4)$ near zero and
	$\Phi(-1/2)=\sqrt3/2+\pi/12-1$. These facts follow from
	$\Phi'(\rho)=\arcsin\rho$, $\Phi''(\rho)=(1-\rho^2)^{-1/2}>0$, and
	$\Phi(0)=0$. Thus any nonzero limiting Gaussian edge correlation gives a
	positive contribution to the covariance of absolute Gaussian edges.

\section{Uniform reduction, profile families, and calibration}
	\label{sec:profiles-calibration}

\subsection{Check of the uniform limit}
	\label{ssec:uniform-limit-check}

	We use the uniform constants
	\[
		A:=1-\frac2\pi,\qquad
		B:=\frac{\sqrt3-2}{\pi}+\frac16
		=\frac2\pi\,\Phi\!\left(-\frac12\right).
	\]

\begin{proposition}\label{prop:uniform-reduction}
	For $p_j \equiv 1/n$ one has $\E D_j=0$, so $d_{\mathbf p}=d$, the uniform statistic.
	Formula (\ref{eq:Var-TVC}) then reduces to the known expression
	\[
		\Var(d_{\mathbf p}) = N\cdot v(n) + O(1)
		= \frac{2N(n-1)}{n}A + \frac{4N(n-2)}{n}B + O(1),
	\]
	where $v(n):=\lim_N\Var(d)/N$ is the uniform leading coefficient of
	\cite[Theorem~3.4]{banic_elezovic_p2}.
\end{proposition}

\noindent
Via Corollary~\ref{cor:n34} the re-derivation covers all $n\ge3$; the uniform theorem
itself covers all $n\ge2$.

\begin{proof}
	Under the uniform $p_j \equiv 1/n$:
	\begin{itemize}
		\item $\Var(D_j) = 2N/n$ for all $j$ (from (\ref{eq:Var-Dj}) with $\delta_j = 0$);
		\item $\rho_j = -1/2$ for all adjacent pairs (from \S\ref{ssec:edge-vars});
		\item $\rho_{j,k} = 0$ for all $|k-j| \ge 2$ (from (\ref{eq:Cov-Dj-disj})),
		and $\Phi(0) = 0$.
	\end{itemize}
	First summand: $(1-2/\pi)(n-1)(2N/n) = \frac{2N(n-1)}{n}A$.
	Adjacent part: $\frac{4}{\pi}(n-2)(2N/n)\Phi(-1/2) = \frac{4N(n-2)}{n}B$.
	The leading disjoint part is zero; its non-Gaussian $O(1)$ residual is contained
	in the remainder of \eqref{eq:Var-TVC}. Summing the displayed leading terms
	matches \cite[Theorem~3.4]{banic_elezovic_p2}.
\end{proof}

	Under the uniform law the edge vector has exactly the covariance of the successive
	differences of $n$ i.i.d.\ $N(0,N/n)$ variables, so the uniform coefficient is the
	classical normal-theory variance of the sum of absolute successive differences
	(Hoel~\cite{hoel1946}, Kamat~\cite{kamat1953}; see \cite[Theorem~3.4]{banic_elezovic_p2},
	where this is recorded): this is why the constants $A$ and $B$ appear, and why the
	large-$n$ constant $2A+4B$ of Proposition~\ref{prop:smooth-universal} is Hoel--Kamat's.

\subsection{Special case: the symmetric comb}
	\label{ssec:comb-closed-form}

	Here the comb is the law under which the test is calibrated: $\mathbf r=\mathbf p$ is the
	comb and $X\sim\Mult(N;\mathbf p)$. For even $n$, consider the symmetric comb (equal
	deviations $\pm\varepsilon$ about $1/n$)
	$p_j = 1/n + (-1)^j \varepsilon$,
	$j = 1, \ldots, n$, $\varepsilon \in [0, 1/n)$ the amplitude of deviation. The range
	$0\le\varepsilon<1/n$ is exactly the one on which all $p_j>0$, as Theorem~\ref{thm:Var-TVC}
	requires; at $\varepsilon=1/n$ a bin vanishes and $\rho_A=-1$. The parity assumption is part
	of the normalization: the alternating perturbations sum to zero only for even
	$n$.

\begin{proposition}[Var($d_{\mathbf p}$) for the symmetric comb]
	\label{prop:comb-closed}
	For the symmetric comb with fixed even $n\ge4$ and fixed
	$0\le\varepsilon<1/n$, put
	\[
		\tilde s^2:=\frac2n-4\varepsilon^2,\qquad
		\rho_A:=\frac{-1/n-\varepsilon+4\varepsilon^2}{\tilde s^2},\qquad
		\rho_B:=\frac{-1/n+\varepsilon+4\varepsilon^2}{\tilde s^2},
	\]
	and
	\[
		\kappa:=\frac{4\varepsilon^2}{\tilde s^2}.
	\]
	Then $\Var(D_j)=N\tilde s^2$ for all $j$; the adjacent correlations alternate between
	$\rho_A$ and $\rho_B$; and for disjoint pairs $|\rho_{j,k}|=\kappa$, with sign
	alternating according to the parity of $k-j$. As $N\to\infty$,
	\begin{multline}\label{eq:Var-TVC-comb}
		\Var(d_{\mathbf p}) = (1 - 2/\pi)(n-1)N\tilde s^2 \\
		+ \frac4\pi N\tilde s^2
		\bigl[\tfrac{n-2}{2}(\Phi(\rho_A) + \Phi(\rho_B))
		+ \tfrac{(n-2)(n-3)}{2}\Phi(\kappa)\bigr] + O(1).
	\end{multline}
\end{proposition}

\noindent
For $\varepsilon\to0$ the formula reduces to Proposition~\ref{prop:uniform-reduction}.
The $O(1)$ remainder in \eqref{eq:Var-TVC-comb} is not uniform in $\varepsilon$
near $1/n$, where $\rho_A\downarrow-1$.

\begin{proof}
	Put $u=1/n$ and
	\[
		v:=\frac{\Var(D_j)}N=2u-4\varepsilon^2.
	\]
	Indeed, every neighboring pair contains one mass $u-\varepsilon$ and one
	mass $u+\varepsilon$, so $p_j+p_{j+1}=2u$, while
	$\delta_j=p_{j+1}-p_j$ alternates between $2\varepsilon$ and
	$-2\varepsilon$. Formula \eqref{eq:Var-Dj} then gives the first assertion.

	For an adjacent pair, substitute the two possible triples
	\[
		(u-\varepsilon,u+\varepsilon,u-\varepsilon)
		\quad\hbox{and}\quad
		(u+\varepsilon,u-\varepsilon,u+\varepsilon)
	\]
	into \eqref{eq:Cov-Dj-adj}. After division by
	$\sigma_j\sigma_{j+1}=N\tilde s^2$, the resulting correlations are
	\[
		\rho_A=\frac{-u-\varepsilon+4\varepsilon^2}{v},
		\qquad
		\rho_B=\frac{-u+\varepsilon+4\varepsilon^2}{v}.
	\]
	Among the $n-2$ adjacent edge pairs, each type occurs $(n-2)/2$ times.

	For disjoint edges, \eqref{eq:Cov-Dj-disj} gives
	\[
		\rho_{j,k}
		=-\frac{\delta_j\delta_k}{v}.
	\]
	Its sign depends on the two edge parities, but its magnitude is the common
	value $\kappa=4\varepsilon^2/v$. Since $\Phi$ is even, only this magnitude
	enters the variance. There are
	\[
		\binom{n-1}{2}-(n-2)=\frac{(n-2)(n-3)}2
	\]
	disjoint pairs. Substitution of these counts and of
	$\sigma_j\sigma_k=N\tilde s^2$ into Theorem~\ref{thm:Var-TVC} gives
	\eqref{eq:Var-TVC-comb}. For $n\ge6$ this is Theorem~\ref{thm:Var-TVC}; for $n=4$ (one
	disjoint pair, the boundary configuration of Remark~\ref{rem:pair-boundary}) it is
	Corollary~\ref{cor:n34}.
\end{proof}

	The formula is even in $\varepsilon$ ($\varepsilon\mapsto-\varepsilon$ swaps $\rho_A$ and
	$\rho_B$), so it also covers the opposite phase $-1/n<\varepsilon<0$: the null variance
	does not depend on which bins are raised.

\begin{remark}[Scaling of the disjoint term]\label{rem:comb-leading}
	For every comb covered by Proposition~\ref{prop:comb-closed},
	$0\le\varepsilon<1/n$, one has
	$\kappa\le2/(n-2)$. Hence $\Phi(\kappa)=O(1/n^2)$ uniformly as $n\to\infty$,
	and the total disjoint contribution is $O(N/n)$, sub-leading relative to the
	per-edge and adjacent terms, which are $O(N)$.
	A concrete example: $n = 6$, a comb with $\varepsilon = 0.136$:
	$\kappa = 0.285$, $\Phi(\kappa) \approx 0.041$. The contributions to $\Var(d_{\mathbf p})/N$:
	per-edge $0.471$, adjacent $0.289$, disjoint $0.081$.
\end{remark}

	\paragraph{Exact validation}
	Proposition~\ref{prop:comb-closed} agrees with exact finite-$N$ variances, computed from the
	exact multinomial marginals of the at most four bins of each edge pair: for $n=6$,
	$\Var(d_{\mathbf p})-Nv(\mathbf p)$ stays bounded for $N=25,50,100$ at
	$\varepsilon=0,\ 1/30,\ 0.136,\ 0.16$ (e.g.\ $-1.07,-1.07,-1.07$ at $\varepsilon=1/(5n)$ and
	$0.77,0.76,0.75$ at $\varepsilon=0.16$), and for $n=4$, $\varepsilon=0.15$ it converges to $0.111$ by
	$N=100$. For the mild comb $\varepsilon=1/(5n)$ one has $\kappa=4/(50n-4)$ and the disjoint part is
	below $0.03\%$ of $v(\mathbf p)$ for every even $n\ge6$; for the strong comb $n=6$,
	$\varepsilon=0.136$ it is $9.6\%$ of $v(\mathbf p)$, and omitting it underestimates
	$\sigma(d_{\mathbf p})/\sqrt N$ by $4.9\%$.

\subsection{Smooth profiles: the Riemann-sum limit}
	\label{ssec:smooth-profile}

	The comb is a periodic reference law. Smooth profiles are obtained by sampling a
	fixed density shape on a finer and finer grid,
	$p_k = f(k/n)/Z_n$, $k=1,\dots,n$, with $f$ a fixed positive $C^1$ function on $[0,1]$ and
	$Z_n=\sum_{k=1}^{n} f(k/n)=n\int_0^1 f+O(1)$ (normalize $\int_0^1 f=1$); the right-endpoint
	sampling $k/n$ is a fixed convention, and any other $O(1/n)$-shifted grid gives the same
	leading coefficient. Linear,
	exponential and other monotone reference laws (ramps, geometric run lengths) are of this type. Here the gradients
	are small:
	$\delta_j=p_{j+1}-p_j=f'(x_j)/n^2+o(n^{-2})=O(n^{-2})$ uniformly at
	$x_j=j/n$, so every edge is locally close to the uniform case and the leading
	variance reduces to Riemann sums.

\begin{proposition}[Smooth-profile universality]\label{prop:smooth-universal}
	Let $f\in C^1([0,1])$ be positive with $\int_0^1 f=1$, and let
	$v_f(n):=\lim_{N\to\infty}\Var(d_{\mathbf p})/N$ be the leading variance coefficient of
	Theorem~\ref{thm:Var-TVC} for $p_k=f(k/n)/Z_n$. Then, as $n\to\infty$,
	\begin{align*}
		\rho_j\to-\tfrac12\ \text{(uniformly)},\qquad
		\rho_{j,k}\to0\ \text{uniformly for } |k-j|\ge2,\\
		v_f(n)\to\sigma_\infty^2=2A+4B,
	\end{align*}
	\emph{independently of the shape} $f$, where $A=1-\tfrac2\pi$ and
	$B=\tfrac{\sqrt3-2}{\pi}+\tfrac16$ (the constants of
	\cite[Definition~3.3]{banic_elezovic_p2}).
\end{proposition}

\begin{proof}
	With $x_j=j/n$ and $p_k=f(x_k)/Z_n=f(x_k)/n+O(n^{-2})$, the adjacent sum and gradient
	are $s_j=p_j+p_{j+1}=2f(x_j)/n+O(n^{-2})$ and $\delta_j=O(n^{-2})$, so
	$\sigma_j^2/N=s_j-\delta_j^2=2f(x_j)/n+O(n^{-2})$ and the per-edge term is a Riemann sum,
	\[
		\Bigl(1-\tfrac2\pi\Bigr)\sum_j\frac{\sigma_j^2}{N}
		=2A\cdot\frac1n\sum_j f(x_j)+O(n^{-1})\longrightarrow 2A\int_0^1 f=2A.
	\]
	For an adjacent pair, the product and quadratic terms of
	\eqref{eq:Cov-Dj-adj} are $O(n^{-2})$. Therefore, uniformly in $j$,
	\[
		\frac{\Cov(D_j,D_{j+1})}{N}
		=-p_{j+1}+O(n^{-2})
		=-\frac{f(x_j)}n+O(n^{-2}),
	\]
	\[
		\frac{\sigma_j\sigma_{j+1}}N
		=\frac{2f(x_j)}n+O(n^{-2}).
	\]
	Hence $\rho_j\to-\tfrac12$ uniformly and
	\[
		\frac4\pi\sum_j\frac{\sigma_j\sigma_{j+1}}{N}\,\Phi(\rho_j)
		\longrightarrow\frac8\pi\,\Phi\bigl(-\tfrac12\bigr)\int_0^1 f=\frac8\pi\,\Phi\bigl(-\tfrac12\bigr)=4B,
	\]
	since $\Phi(-\tfrac12)=\tfrac{\sqrt3}2+\tfrac{\pi}{12}-1$ gives
	$\tfrac8\pi\Phi(-\tfrac12)=\tfrac{4\sqrt3-8}{\pi}+\tfrac23=4B$.
	Finally, for disjoint edges,
	\[
		\rho_{j,k}
		=-\frac{\delta_j\delta_k}
		{\sqrt{(s_j-\delta_j^2)(s_k-\delta_k^2)}}=O(n^{-3})
	\]
	uniformly, because positivity of $f$ makes the denominator of order $n^{-1}$.
	Thus $\Phi(\rho_{j,k})=O(n^{-6})$, while
	$\sigma_j\sigma_k/N=O(n^{-1})$. There are $O(n^2)$ disjoint pairs, so their
	total contribution to $v_f(n)$ is $O(n^{-5})$. Summing the three parts gives
	$v_f(n)\to2A+4B=\sigma_\infty^2$.
\end{proof}

\begin{remark}[Numerical next-order behavior]\label{rem:smooth-1n}
	The limit is the \emph{same} universal constant $\sigma_\infty^2$ as the uniform
	large-$n$ value $v(n)=\sigma_\infty^2-C_1/n$, where
	$C_1=2A+8B$ is the linear--cyclic gap constant of
	\cite[Definition~3.3 and Theorem~3.4]{banic_elezovic_p2}: the
	density shape does not enter the leading term, because only $\int_0^1 f=1$ enters the
	two Riemann sums. Numerically, in the tested profile families (linear ramps of
	several slopes; exponential and quadratic profiles) the next order is
	consistent with
	$v_f(n)=\sigma_\infty^2+C_1[f]/n+O(n^{-2})$. Every tested \emph{linear} ramp gives
	$C_1[f]=-C_1=-(2A+8B)$ to all tested digits, whereas the numerics for curved
	profiles (exponential, quadratic) are consistent with an additional shape
	term; the experiments rule out proportionality to the Fisher functional
	$\int(f')^2/f$, but we do not conjecture a closed form here.
\end{remark}

\subsection{Application: test calibration}
	\label{ssec:test-calibration-nonuniform}

	This subsection is an application, downstream of the probabilistic content above
	(Corollary~\ref{cor:weak-limit}): formula (\ref{eq:Var-TVC}) supplies the second input to a practical
	approximation of the distribution of $d_{\mathbf p}$ under a non-uniform $H_0$. The
	recipe is the moment-matched Gamma calibration introduced for the comb test
	in~\cite[Sec.~IV-C]{banic_elezovic_ct}, a \emph{two-moment Gamma surrogate}: fit a Gamma law by matching the
	mean --- from \cite[Theorem~4.3]{banic_elezovic_p3}, including its oscillatory (lattice) term ---
	and the variance from \eqref{eq:Var-TVC}, and reject in the upper tail. For rational $\mathbf r$ and
	moderate $N$ the exact null law of $d_{\mathbf r}$ is computable by a weighted form of the dynamic
	programme of~\cite[Eq.~(4)]{banic_elezovic_ct}; the surrogate is meant for the regime where that is too costly. We
	emphasize that the moment calculation does not justify a Gamma \emph{shape} for
	$d_{\mathbf p}$, whose fixed-$n$ limit law is a functional of $n-1$ correlated folded
	Gaussians; the fit is a moment-matched heuristic in the classical tradition of
	two-moment chi-square/Gamma approximations (Patnaik~\cite{patnaik1949}), to be
	justified or replaced by a dedicated distributional analysis. The calibrations
	below are illustrations resting on the leading fixed-$n$ coefficient
	$v(\mathbf p)$ of Theorem~\ref{thm:Var-TVC}; the large-$n$ entries ($n=100$)
	invoke the iterated-limit reading of Remark~\ref{rem:VarTVC-scope} (first
	$N\to\infty$, then $n\to\infty$), not a joint uniform-in-$n$ error bound.

	\paragraph{Role of the disjoint terms}
	The disjoint covariances must not be dropped for strongly non-uniform nulls: in the example of
	Remark~\ref{rem:comb-leading} they are about a tenth of the leading variance (a fifth of its covariance
	part), and omitting them understates $\sigma(d_{\mathbf p})$ by $5\%$. For the moderate comb of
	Table~\ref{tab:typeI} they are about $1\%$ of the variance and do not move the critical value; there the
	table shows the effect of the adjacent covariances, which raise the size from $5\%$ to $8$--$9\%$ when
	omitted. Near-uniform caveat: when $\mathbf p$ is close to uniform the leading disjoint term is tiny
	($\Phi(\rho)\approx\rho^2/2$) and is dominated by the $O(1)$ disjoint remainders, which in the uniform case
	sum to about $-n/\pi$ (Part~2), a relative correction of order $1/\bar m$; at small $\bar m$ the
	second-order uniform constants of~\cite{banic_elezovic_p2} are the better approximation to
	$\Var(d_{\mathbf p})$ --- the statistic remains the centered $d_{\mathbf p}$ --- and for rational
	$\mathbf r$ the exact null law is available by a weighted form of the dynamic programme of
	\cite[Eq.~(4)]{banic_elezovic_ct} (see \S\ref{ssec:test-calibration-nonuniform} below).

	\paragraph{Empirical Type I error} For a moderate comb
	($p_{\text{low}} = 0.5/n$, $p_{\text{high}} = 1.5/n$, a $\pm50\%$ modulation; the mild comb of
	\S\ref{ssec:comb-closed-form} has $\varepsilon=1/(5n)$), nominal $\alpha = 5\%$,
	upper-tail test; $\bar m:=N/n$ is the expected count per bin on the uniform
	scale.
	\begin{table}[ht]
	\centering
	\footnotesize
	\setlength{\tabcolsep}{4pt}%
	\begin{tabular}{l|c|c|c|c|c}
		$(n,\bar m)$ & Gamma, per-edge & Gamma, no disj.\ & Gamma, full & normal & $\chi^2_{n-1}$ \\\hline
		$(10,30)$    & $8.30\%$ & $4.69\%$ & $4.69\%$ & $6.29\%$ & $5.00\%$ \\
		$(10,100)$   & $8.28\%$ & $4.97\%$ & $4.97\%$ & $6.13\%$ & $5.00\%$ \\
		$(100,300)$  & $9.27\%$ & $4.93\%$ & $4.88\%$ & $5.42\%$ & $5.03\%$ \\
		$(100,1000)$ & $9.23\%$ & $4.89\%$ & $4.87\%$ & $5.39\%$ & $5.04\%$ \\
	\end{tabular}
	\caption{Realized size at nominal $5\%$ under the moderate-comb null ($p_{\mathrm{low}}=0.5/n$,
	$p_{\mathrm{high}}=1.5/n$), upper-tail test, $\bar m=N/n$. Columns 2--4: two-moment Gamma surrogate with
	the mean of \cite[Theorem~4.3]{banic_elezovic_p3} and, respectively, the per-edge variance
	$(1-\tfrac2\pi)\sum_j\Var(D_j)$, the variance \eqref{eq:Var-TVC} without its disjoint pairs, and the full
	\eqref{eq:Var-TVC}; column 5: normal cut with the same two moments; column 6: Pearson's $\chi^2$ with its
	asymptotic $\chi^2_{n-1}$ calibration~\cite{pearson1900}. Monte Carlo with $K=4\cdot10^5$ replications for
	$n=10$ and $2\cdot10^5$ for $n=100$ (binomial standard error $0.03\%$, resp.\ $0.05\%$).}
	\label{tab:typeI}
	\end{table}
	With \eqref{eq:Var-TVC} the two-moment Gamma surrogate has size $4.7$--$5.0\%$ at nominal $5\%$ in every
	tested configuration; it is slightly conservative, and more so at smaller levels (size $0.8$--$0.95\%$ at
	nominal $1\%$, $0.06$--$0.09\%$ at $0.1\%$), because the moment-matched Gamma has a heavier right tail than
	the limit law $W$ of Corollary~\ref{cor:weak-limit} (skewness $2\sqrt{v}/\mu\approx0.62$ against $0.51$ for
	$W$ at $n=10$; $0.19$ against $0.16$ at $n=100$). A normal cut with the same two moments is
	anti-conservative ($5.4$--$6.3\%$ at nominal $5\%$, $1.3$--$2.0\%$ at $1\%$): a reader holding a mean and a
	variance should not reach for $z$-scores. Dropping all covariance terms (per-edge variance only) inflates
	the size to $8$--$9\%$; on this comb the disjoint pairs are $1\%$ ($n=10$) and $0.14\%$ ($n=100$) of the
	variance and leave the size unchanged, so what the table displays is the role of the \emph{adjacent}
	covariances. Pearson's $\chi^2$ is at nominal in all rows. The fixed-$n$ limit law is a sum of correlated
	folded Gaussians, not a Gamma law; the surrogate is a moment-matched approximation whose accuracy is
	documented here, not derived.

\section{Variance of the reference-centered statistic $d_{\mathbf r}$}
	\label{sec:Var-dXY}

	This section computes the variance of the reference-centered statistic $d_{\mathbf r}$.
	It belongs here, rather than in a sequel, because it is the same folded-normal
	covariance calculation seen from the opposite end. In the matched statistic
	$d_{\mathbf p}$ the fold of every edge is placed at its mean, so the absolute
	values do not linearize and the folded-normal covariance terms remain. For a
	fixed reference $\mathbf r$ different from the sampling law $\mathbf p$, an
	edge with $\Delta_i\ne0$ has its fold displaced from the mean; then the
	absolute value linearizes on that edge, and the bivariate folded covariance of
	Proposition~\ref{prop:bivar-folded} degenerates to an ordinary multinomial quadratic
	form. The section thus provides the bridge between the matched
	$d_{\mathbf p}$ regime and fixed alternatives, and it is where the general
	bivariate folded-normal covariance, a classical quantity
	(Kan--Robotti~\cite{kanrobotti2017}), is derived in the integral form used below.

	Let the observed counts follow $X\sim\Mult(N;\mathbf p)$ and let $\mathbf r$ be a fixed
	reference distribution in the open simplex: $r_i>0$ for all $i$ and $\sum_ir_i=1$ (for the
	triangular-array result the sequence $\mathbf r^{(N)}$ likewise lies in the simplex).
	The boundary faces $p_i=0$ are excluded because they require a separate
	degenerate edge analysis; the restriction $r_i>0$ keeps the reference law in
	the same open-simplex regime as the matched theory. For edge $i$ define the gradient
	mismatch
	\[
		\Delta_i:=(p_{i+1}-p_i)-(r_{i+1}-r_i)
	\]
	(denoted $\delta_i$ in \cite{banic_elezovic_p3}; here the letter $\delta$ is
	reserved for the gradient $\delta_j=p_{j+1}-p_j$ of \S\S2--4) and the
	reference-centered edge
	\[
		W_i:=D_i-N(r_{i+1}-r_i)=\tilde D_i+N\Delta_i,\qquad
		\tilde D_i:=D_i-N(p_{i+1}-p_i).
	\]
	The reference-centered statistic is
	\[
		d_{\mathbf r}:=\sum_{i=1}^{n-1}|W_i|
	\]
	This is the statistic whose expectation was studied
	in~\cite[\S6]{banic_elezovic_p3}. Here we compute its variance; throughout this
	section
	$\mu_i:=N\Delta_i$ and $\sigma_i^2:=\Var(D_i)$.
	A nonzero mean inside the absolute values \emph{simplifies} the
	variance (and, when it is the uncentered statistic on a non-uniform null, this
	simplification is the loss of phase-blindness noted in the introduction): the
	absolute value linearizes, so $\Var(d_{\mathbf r})$ acquires a
	closed form considerably simpler than $\Var(d_{\mathbf p})$.

	Each edge falls into one of three regimes, according to the size of the standardized
	mismatch $a_i=N\Delta_i/\sigma_i$, and the results below are organized by them:
	\[
		\underbrace{\Delta_i=0}_{\text{matched}}
		\qquad
		\underbrace{a_i\to a_i^*\in\R}_{\text{local (array $\mathbf r^{(N)}$)}}
		\qquad
		\underbrace{\Delta_i\ne0\text{ fixed}\ \Rightarrow\ |a_i|\to\infty}_{\text{fixed mismatch}} .
	\]
	The matched and fixed-mismatch regimes are those of a fixed pair $(\mathbf p,\mathbf r)$;
	Theorem~\ref{thm:Var-dXY} covers a fully fixed mismatch and carries an exponentially small
	error, while the local regime is genuinely a triangular-array statement and is treated,
	with an $o(N)$ error, in Proposition~\ref{prop:dXY-mixed}. A single edge may sit in any
	one of the three; the theorems below allow the regimes to be mixed across edges.

	\paragraph{Per-edge variance: the folded-normal profile $V$} Since
	$|W_i|^2 = W_i^2$, the second moment is exact: $\E|W_i|^2 = \mu_i^2 + \sigma_i^2$.

\begin{lemma}[Per-edge folded variance]\label{lem:dXY-Vprofile}
	Define
	\[
		V(a):=1+a^2-G(a)^2,\qquad a\in\R,
	\]
	where $G$ is the folded-normal profile of
	Definition~\ref{def:folded-profile}. Let
	$W_i^{(G)}\sim\mathcal N(\mu_i,\sigma_i^2)$ and put
	$a_i:=\mu_i/\sigma_i=N\Delta_i/\sigma_i$. Then
	\begin{equation}\label{eq:dXY-peredge-var}
		\Var(|W_i^{(G)}|)=\sigma_i^2\,V(a_i).
	\end{equation}
	For the multinomial edge variable $W_i=\tilde D_i+N\Delta_i$, at fixed
	$\mathbf p,\mathbf r$ and fixed $n$ with all $p_i>0$,
	\[
		\Var(|W_i|)=\sigma_i^2V(a_i)+o(N),\qquad N\to\infty .
	\]
	Moreover, $V$ is even, $V(0)=1-2/\pi$, $V$ is strictly increasing on
	$[0,\infty)$ from $1-2/\pi$ to $1$, and
	\[
		1-V(a)=O(|a|^{-1}e^{-a^2/2}),\qquad |a|\to\infty ,
	\]
	more precisely $1-V(a)\sim4\phi(a)/a$.
\end{lemma}

\noindent
For fixed $\mathbf p,\mathbf r$, each edge is either matched ($\Delta_i=0$) or
has a fixed mismatch ($\Delta_i\ne0$). The local regime, in which $a_i$ has a
finite nonzero limit, requires a triangular array
$\mathbf r=\mathbf r^{(N)}$ and is treated in
Proposition~\ref{prop:dXY-mixed}. The profile $G$ is the folded-normal mean
profile used in \cite[Definition~6.3, Lemma~6.4]{banic_elezovic_p3} and goes
back to Leone et al.~\cite{leone1961}; the value $V(0)=1-2/\pi$ is the constant
$A$ of \cite[Definition~3.3]{banic_elezovic_p2}.

\begin{proof}
	The displayed identity follows from
	$\Var|W_i^{(G)}| = (\mu_i^2 + \sigma_i^2) - \sigma_i^2 G(a_i)^2$.
	For a matched edge, $\Delta_i=0$ and $a_i=0$. The standardized edge
	$\tilde D_i/\sigma_i$ converges to $Z\sim\mathcal N(0,1)$, and bounded
	multinomial increments give uniform integrability of its square. Hence the
	first two absolute moments converge and
	$\Var(|W_i|)=\sigma_i^2V(0)+o(N)$.

	For a fixed mismatch, $\Delta_i\ne0$ and therefore
	$|a_i|=N|\Delta_i|/\sigma_i\asymp\sqrt N$. The large-deviation
	linearization used in Theorem~\ref{thm:Var-dXY} gives
	$\Var(|W_i|)=\sigma_i^2+o(N)$. On the other hand,
	$V(a_i)\to1$, so also $\sigma_i^2V(a_i)=\sigma_i^2+o(N)$.
	These are the only two possibilities for fixed $\mathbf p,\mathbf r$.
	Evenness from the folded-profile properties
	\cite[Lemma~6.4]{banic_elezovic_p3}; $V(0) = 1 - G(0)^2 = 1 - 2/\pi$; from
	$G(a) = |a| + O(e^{-a^2/2})$ it follows that $V(a) \to 1$, and the sharp rate
	follows from the Mills expansion $G(a)-|a|\sim2\phi(a)/a^2$ together with
	$V(a)=1-(G(a)-|a|)(G(a)+|a|)$.
	For monotonicity, let $W_a=Z+a$ with $Z\sim\mathcal N(0,1)$. Since
	$G'(a)=\E[\operatorname{sgn}(W_a)]$ --- differentiation under the expectation
	being justified by dominated convergence applied to the a.e.-derivative of the
	$1$-Lipschitz map $a\mapsto|z+a|$ --- and
	$\E[|W_a|\operatorname{sgn}(W_a)]=\E W_a=a$,
	\[
		\frac12V'(a)=a-G(a)G'(a)
		=\Cov\bigl(|W_a|,\operatorname{sgn}(W_a)\bigr).
	\]
	For $a>0$, the conditional density of $W_a$ on $(0,\infty)$ and that of
	$-W_a$ conditional on $W_a<0$ have likelihood ratio proportional to
	$e^{2ax}$, which is strictly increasing in $x>0$. Hence
	\[
		m_+:=\E(W_a\mid W_a>0)>
		m_-:=\E(-W_a\mid W_a<0).
	\]
	If $p_\pm=\Pr(\pm W_a>0)$, direct conditioning on the sign gives
	\[
		\Cov\bigl(|W_a|,\operatorname{sgn}(W_a)\bigr)
		=2p_+p_-(m_+-m_-)>0.
	\]
	Therefore $V'(a)>0$ for $a>0$.
\end{proof}

\begin{remark}[Interpretation of $V$]\label{rem:V-interp}
	$V$ measures the effect of folding on the variance. On a matched edge
	($a=0$) the variance is reduced to $A\approx0.363$ of $\sigma_i^2$
	(half-normal case). On an edge with a strong mismatch ($|a|\to\infty$),
	$W_i$ has a fixed sign with probability tending to one, so $|W_i|=\pm W_i$
	apart from an exponentially small crossing event and $V(a)\to1$.
\end{remark}

	\paragraph{Covariances: bivariate folded normal} For a pair $(W_i,W_j)$,
	the limiting normal law has means $(\mu_i,\mu_j)$, variances
	$(\sigma_i^2,\sigma_j^2)$, and correlation
	\[
		\rho_{ij}=\frac{\Cov(D_i,D_j)}{\sigma_i\sigma_j}.
	\]
	Centering by constants does not change this covariance. The general
	$\Cov(|W_i|,|W_j|)$ is a bivariate folded normal, a function of
	$(a_i,a_j,\rho_{ij})$ expressible through the bivariate normal CDF
	(Kan--Robotti~\cite{kanrobotti2017}).
	Two limiting regimes are used below:
	\begin{itemize}
		\item \emph{matched pair} ($a_i = a_j = 0$): Nabeya
		(Corollary~\ref{cor:cov-abs}),
	\[
		\Cov(|W_i|, |W_j|) =
		\tfrac{2}{\pi}\sigma_i\sigma_j\Phi(\rho_{ij})+o(N);
	\]
		the sum reproduces the leading term of $\Var(d_{\mathbf p})$
		(Theorem~\ref{thm:Var-TVC}).
		\item \emph{strongly-mismatched pair} ($|a_i|, |a_j| \to \infty$): linearization
		$|W_i| = |\mu_i| + \epsilon_i\tilde D_i + (\text{exp.\ small})$,
		$\epsilon_i := \operatorname{sgn}(\Delta_i)$, so
	\begin{equation}\label{eq:dXY-cov-mismatch}
		\Cov(|W_i|, |W_j|) =
		\epsilon_i\epsilon_j\,\Cov(D_i, D_j)+O\bigl(e^{-c'N}\bigr)
	\end{equation}
		--- the absolute value linearizes, and no Nabeya $\Phi$ appears
		(the exponential error comes from the linearization bound in the proof of
		Theorem~\ref{thm:Var-dXY}).
	\end{itemize}

	The two regimes are endpoints of a single closed form obtained from Price's
	theorem.

\begin{proposition}[Bivariate folded-normal covariance]\label{prop:bivar-folded}
	Let $(W_i,W_j)$ be jointly normal with means $\mu_i=\sigma_i a_i$,
	$\mu_j=\sigma_j a_j$, variances $\sigma_i^2,\sigma_j^2$ and correlation
	$|\rho|<1$. Let $\Phi_{\mathcal N}$ be the standard normal CDF, and let
	$\Phi_{\mathcal N_2}(\cdot,\cdot;\rho)$ denote the standard bivariate normal
	CDF with correlation $\rho$, with the convention
	\[
		\Pr(Z_i+a_i>0,\,Z_j+a_j>0)=\Phi_{\mathcal N_2}(a_i,a_j;\rho).
	\]
	Then
	\begin{equation}\label{eq:bivar-folded}
		\Cov(|W_i|,|W_j|)
		= \sigma_i\sigma_j\Bigl[\,4\!\int_0^{\rho}\Phi_{\mathcal N_2}(a_i,a_j;t)\,dt
		- \rho\bigl(2\Phi_{\mathcal N}(a_i)+2\Phi_{\mathcal N}(a_j)-1\bigr)\Bigr],
	\end{equation}
	where the integral is oriented when $\rho<0$.
\end{proposition}

\noindent
The notation $\Phi_{\mathcal N_2}$ is kept distinct from the Nabeya
function $\Phi$ of \eqref{eq:SM-Phi}. The integral in
\eqref{eq:bivar-folded} can also be written in terms of
$\Phi_{\mathcal N}$ and $\phi$ using the folded-moment
formulas of Kan--Robotti~\cite{kanrobotti2017}. In applications below it is
evaluated at the edge correlation $\rho=\rho_{ij}$.

\begin{proof}
	Write $W_i=\sigma_i(Z_i+a_i)$ with $(Z_i,Z_j)$ standard bivariate normal of
	correlation $\rho$, so $\Cov(|W_i|,|W_j|)=\sigma_i\sigma_j\,[M(\rho)-G(a_i)G(a_j)]$
	with $M(\rho):=\E\bigl[|Z_i+a_i|\,|Z_j+a_j|\bigr]$. The mixed distributional
	derivative
	$\partial_U\partial_V\,|U|\,|V|=\operatorname{sgn}(U)\operatorname{sgn}(V)$ carries no
	atom or singular part --- the kinks of $|U|$ and $|V|$ sit on separate axes,
	and the mixed derivative is a bounded measurable function --- so the classical
	form of Price's theorem~\cite{price1958} applies:
	\[
		\frac{dM}{d\rho}=\E\bigl[\operatorname{sgn}(Z_i+a_i)\operatorname{sgn}(Z_j+a_j)\bigr]
		=4\,\Phi_{\mathcal N_2}(a_i,a_j;\rho)-2\Phi_{\mathcal N}(a_i)-2\Phi_{\mathcal N}(a_j)+1,
	\]
	using the identity $\Pr(Z_i+a_i>0,Z_j+a_j>0)=\Pr(-Z_i<a_i,-Z_j<a_j)
		=\Phi_{\mathcal N_2}(a_i,a_j;\rho)$ (valid because $(-Z_i,-Z_j)$ has the same
		correlation $\rho$, so the upper-tail event equals the lower-tail bivariate
		normal CDF at $(a_i,a_j)$) and
	$\operatorname{sgn}=2\,\mathbb 1_{>0}-1$. At $\rho=0$ independence gives
	$M(0)=G(a_i)G(a_j)$; integrating from $0$ to $\rho$ and subtracting $G(a_i)G(a_j)$
	gives \eqref{eq:bivar-folded}.
\end{proof}

\begin{remark}[Reduction to the two regimes]\label{rem:bivar-limits}
	Formula \eqref{eq:bivar-folded} interpolates the cases already on record.
	\emph{Matched} ($a_i=a_j=0$): $\Phi_{\mathcal N_2}(0,0;t)=\tfrac14+\tfrac{\arcsin t}{2\pi}$
	and $\Phi_{\mathcal N}(0)=\tfrac12$ reduce it to
	$\tfrac2\pi\sigma_i\sigma_j(\sqrt{1-\rho^2}+\rho\arcsin\rho-1)=\tfrac2\pi\sigma_i\sigma_j\Phi(\rho)$,
	the Nabeya value. \emph{Strongly mismatched} ($|a_i|,|a_j|\to\infty$):
	$\Phi_{\mathcal N_2}(a_i,a_j;t)\to\mathbb 1\{a_i>0\}\mathbb 1\{a_j>0\}$ and
	$\Phi_{\mathcal N}(a_i)\to\mathbb 1\{a_i>0\}$ collapse it to
	$\sigma_i\sigma_j\,\epsilon_i\epsilon_j\,\rho=\epsilon_i\epsilon_j\Cov(D_i,D_j)$
	($\epsilon=\operatorname{sgn}(a)$), recovering \eqref{eq:dXY-cov-mismatch}. Together
	with the per-edge $\sigma_i^2V(a_i)$ (Lemma~\ref{lem:dXY-Vprofile}), this completes
	the \emph{Gaussian-limit} variance functional for an arbitrary mixture of
	matched, local and mismatched edges; the transfer of this functional to the
	multinomial statistic is the content of Proposition~\ref{prop:dXY-mixed}
	below. Numerically
	\eqref{eq:bivar-folded} matches direct integration to $10^{-12}$ across matched,
	mixed, and mismatched pairs.
\end{remark}

	We first record the simplest fixed-alternative case, where every edge has a
	nonzero mismatch and the absolute values linearize with exponentially small
	error.

\begin{theorem}[Variance at a fully mismatched fixed alternative]\label{thm:Var-dXY}
	Let $n\ge2$, and let $\mathbf p, \mathbf r$ be fixed with $\Delta_i \ne 0$ for
	all $i$ (every
	edge strongly mismatched) and $p_i>0$ for all $i$. Then for $N \to \infty$:
	\begin{align}\label{eq:Var-dXY-closed}
		\Var(d_{\mathbf r}) &= \Var\Bigl(\sum_{i=1}^{n-1}\epsilon_i D_i\Bigr)
		+ O\bigl(e^{-c'N}\bigr)\nonumber\\
		&= N\Bigl[\sum_{k=1}^{n} c_k^2\, p_k - \Bigl(\sum_{k=1}^{n} c_k\, p_k\Bigr)^2\Bigr]
		+ O\bigl(e^{-c'N}\bigr),
	\end{align}
	where $\epsilon_i = \operatorname{sgn}(\Delta_i)$ and
	$c_k := \epsilon_{k-1} - \epsilon_k \in \{-2, \ldots, 2\}$
	($\epsilon_0 := \epsilon_n := 0$).
\end{theorem}

	The statement covers every $n\ge2$. The proof below is given for $n\ge3$, the
	scope of the surrounding sections; the two-bin case, in which the statistic
	collapses to a folded binomial $d_{\mathbf r}=2|X_1-Nr_1|$, is settled in
	\S\ref{ssec:small-n2}, in a sharper form that makes the exponentially small
	remainder explicit.

\begin{proof}
	From the linearization,
	\[
		d_{\mathbf r} = \sum_i|\mu_i| + \sum_i\epsilon_i\tilde D_i + (\text{exp.\ small}),
	\]
	so $\Var(d_{\mathbf r}) = \Var(\sum_i\epsilon_i D_i) + O(e^{-c'N})$ (the
	deterministic
	term and the centering do not change the variance; moreover the $(\text{exp.\ small})$
	remainder $R_i := |W_i| - |\mu_i| - \epsilon_i\tilde D_i = 2(\epsilon_i W_i)^- \ge 0$ is supported on the crossing event $\{\epsilon_iW_i<0\}$, on which
	$R_i$ can be as large as $O(N)$ while the event has probability $O(e^{-cN})$
		by Hoeffding's inequality for sums of bounded
	increments (the crossing threshold is the fixed positive multiple
	$N|\Delta_i|$ of $N$); hence $\E[R_i^2]\le CN^2e^{-cN}$, and by Cauchy--Schwarz every
	covariance it contributes is $O(N^{3/2}e^{-cN/2})=O(e^{-c'N})$, negligible in $L^2$). By the expansion
	$\sum_i\epsilon_i(X_{i+1} - X_i) = \sum_k(\epsilon_{k-1} - \epsilon_k)X_k = \sum_k c_k X_k$,
	by the multinomial variance of a linear combination (from (\ref{eq:mult-cov})):
	$\Var(\sum_k c_k X_k) = N[\sum_k c_k^2 p_k - (\sum_k c_k p_k)^2]$.
\end{proof}

	The constant $c'$ degenerates as $\min_i|\Delta_i|\to0$ (Hoeffding on the crossing event
	gives $\Pr\le e^{-N\Delta_i^2/2}$, so any $c'<\min_i\Delta_i^2/4$ works); for moderate $N$,
	or when some standardized mismatch $|a_i|=N|\Delta_i|/\sigma_i$ is below about $3$, the
	functional \eqref{eq:Var-dXY-general} evaluated at the actual $a_i$ is the accurate
	approximation (in an $n=3$ example with $|\Delta_i|=1/30$ the closed form is $16\%$ off at
	$N=640$, the functional $0.09\%$).

\begin{remark}[Why it is simpler than $\Var(d_{\mathbf p})$]\label{rem:dXY-simpler}
	$\Var(d_{\mathbf p})$ required the Nabeya function $\Phi(\rho)$, per-edge
	$1/N$ corrections and disjoint covariances (\S\ref{ssec:Var-TVC-main}). For a
	fully mismatched fixed alternative these folded-normal corrections are replaced by
	linear covariance terms:
	$\Var(d_{\mathbf r})$ is the variance of a \emph{linear} combination of multinomial
	components --- an elementary expression without transcendental constants.
	The reason is the fixed-sign linearization of every absolute value
	(Remark~\ref{rem:V-interp}). The same linearization with $\mathbf r=\mathbf u$ and a
	non-uniform sampling law $\mathbf p$ --- every $\Delta_i=p_{i+1}-p_i\ne0$ --- is the
	sign locking that makes the uncentered statistic $d_{\mathbf u}$ a linear functional
	$\sum_i\epsilon_iD_i$ of the counts and so lose the phase-blindness of the comb test;
	the simplification is precisely what centering at the null avoids
	(\cite[\S6]{banic_elezovic_p3}).
\end{remark}

\begin{proposition}[Mixed regimes: triangular array]\label{prop:dXY-mixed}
	Fix $\mathbf p$ with all $p_i>0$ and fixed $n\ge3$, and let the reference vary
	with $N$, $\mathbf r=\mathbf r^{(N)}$, in such a way that every standardized
	mismatch converges, $a_i^{(N)}=N\Delta_i^{(N)}/\sigma_i\to a_i^*\in[-\infty,+\infty]$.
	Let $\mathcal C(a,b;\rho)$ denote the normalized bivariate folded-normal
	covariance from Proposition~\ref{prop:bivar-folded}, so that
	$\sigma_i\sigma_j\mathcal C(a,b;\rho)$ is the covariance in
	\eqref{eq:bivar-folded}; at infinite arguments $\mathcal C$ is defined by its
	limiting value. Then
	\begin{equation}\label{eq:Var-dXY-general}
		\Var(d_{\mathbf r}) = \sum_{i=1}^{n-1}\sigma_i^2\,V\bigl(a_i^{(N)}\bigr)
		+ 2\!\!\sum_{1 \le i < j \le n-1}\!\!
		\mathcal C\bigl(a_i^{(N)},a_j^{(N)};\rho_{ij}\bigr)\,\sigma_i\sigma_j + o(N),
	\end{equation}
\end{proposition}

\noindent
The limiting values of $\mathcal C$ exist by dominated convergence in
\eqref{eq:bivar-folded}: for one infinite argument,
$\mathcal C(\pm\infty,a;\rho)=\pm\rho\,G'(a)$, and for two infinite arguments
$\mathcal C(\pm\infty,\pm\infty;\rho)=\epsilon_i\epsilon_j\rho$, with signs
chosen according to the two limits. Thus, if all $a_i^*=0$, the display reduces
to the leading term of $\Var(d_{\mathbf p})$; if all $|a_i^*|=\infty$, it reduces to the
leading linearized variance $\Var(\sum_i\epsilon_iD_i)+o(N)$, namely the leading
term of \eqref{eq:Var-dXY-closed} without its exponential remainder. That sharper
error is a separate fixed-alternative result, valid when the $\Delta_i$ stay
bounded away from $0$.

\begin{proof}
	The multinomial increments are bounded, so all moments of the standardized
	edge pairs are bounded uniformly in $N$. Hence, whenever the $a_i^{(N)}$ are
	bounded, the normalized products
	\[
		\frac{|W_iW_j|}{\sigma_i\sigma_j}
		=\left|\frac{\tilde D_i}{\sigma_i}+a_i^{(N)}\right|
		 \left|\frac{\tilde D_j}{\sigma_j}+a_j^{(N)}\right|
	\]
	and the squares $W_i^2/\sigma_i^2$ are uniformly integrable. On edges and
	pairs with bounded standardized mismatch, the bivariate
	central limit theorem with the deterministic drifts $a_i^{(N)}\to a_i^*$
	therefore gives convergence of the first two absolute moments and of the mixed
	absolute moment to their Gaussian values, i.e.\ each per-edge variance is
	$\sigma_i^2V(a_i^{(N)})+o(N)$ and each per-pair covariance is the bivariate
	folded value $+\,o(\sigma_i\sigma_j)=o(N)$. On edges with
	$|a_i^{(N)}|\to\infty$ the fixed-alternative linearization of
	Theorem~\ref{thm:Var-dXY} is \emph{not} available verbatim: for slowly diverging
	$a_i^{(N)}$ the crossing threshold is only
	$N|\Delta_i^{(N)}|=|a_i^{(N)}|\sigma_i$, of order $|a_i^{(N)}|\sqrt N$ rather than
	$N$, so
	the $e^{-cN}$ estimate fails. We use instead a standardized sub-Gaussian tail
	bound. With $Z_i^{(N)}=\tilde D_i/\sigma_i$ and
	$\epsilon_i=\operatorname{sgn}a_i^{(N)}$,
	\[
		\frac{|W_i|}{\sigma_i}-|a_i^{(N)}|
		=\epsilon_iZ_i^{(N)}+2\bigl(-\epsilon_iZ_i^{(N)}-|a_i^{(N)}|\bigr)_+ .
	\]
	Since $\tilde D_i$ is a sum of bounded centered increments with
	$\sigma_i^2\asymp N$, Hoeffding gives
	$\Pr\{\epsilon_iZ_i^{(N)}<-t\}\le Ce^{-ct^2}$ with constants independent of $N$,
	whence
	\[
		\E\bigl[\bigl(-\epsilon_iZ_i^{(N)}-|a_i^{(N)}|\bigr)_+^2\bigr]
		=\int_0^\infty 2u\,\Pr\{\epsilon_iZ_i^{(N)}<-|a_i^{(N)}|-u\}\,du\longrightarrow0 .
	\]
	Thus $|W_i|/\sigma_i-\epsilon_iZ_i^{(N)}\to0$ in $L^2$ (modulo the deterministic
	$|a_i^{(N)}|$), and after multiplying by $\sigma_i\asymp\sqrt N$ every variance and
	covariance error involving such an edge is $o(N)$; this yields the limiting
	entries of Remark~\ref{rem:bivar-limits} at the leading $o(N)$ level, without the
	sharp exponential error of the fixed-alternative theorem. Summing the fixed number of
	terms gives \eqref{eq:Var-dXY-general}.
\end{proof}

	The $o(N)$ cannot in general be improved to $O(1)$: an edge with a nonzero standardized
	mismatch couples to the third cumulants of the counts, and the remainder is of exact
	order $\sqrt N$ --- for $n=3$, $\mathbf p=(\tfrac14,\tfrac14,\tfrac12)$, $\mathbf r=\mathbf u$
	(one matched, one mismatched edge) the exact variance minus \eqref{eq:Var-dXY-general} is
	$-1.5,-2.5,-3.7,-5.3,-7.5,-10.6$ at $N=20,40,\dots,640$, growing like $\sqrt2$ per
	doubling --- in agreement with the $O(1)$ skewness correction to $\E(d_{\mathbf r})$ in
	\cite[Corollary~6.9]{banic_elezovic_p3}.

\begin{corollary}[Fixed pair with matched and mismatched edges]\label{cor:dXY-fixed-mixed}
	Fix $n\ge3$, an interior $\mathbf p$ and a reference $\mathbf r$, and let
	$M:=\{i:\Delta_i=0\}$ be the matched edges. Then
	\begin{multline*}
		\Var(d_{\mathbf r})=\Bigl(1-\tfrac2\pi\Bigr)\sum_{i\in M}\sigma_i^2+\sum_{i\notin M}\sigma_i^2
		+\frac4\pi\sum_{\substack{i<j\\ i,j\in M}}\sigma_i\sigma_j\Phi(\rho_{ij})\\
		+2\sum_{\substack{i<j\\ i,j\notin M}}\epsilon_i\epsilon_j\Cov(D_i,D_j)+o(N):
	\end{multline*}
	matched--mismatched cross covariances vanish at leading order, since
	$\mathcal C(0,\pm\infty;\rho)=\pm\rho\,G'(0)=0$.
\end{corollary}

\begin{proof}
	Proposition~\ref{prop:dXY-mixed} with $a_i^*=0$ on $M$ and $\pm\infty$ off $M$;
	$G'(0)=\operatorname{erf}(0)=0$.
\end{proof}

	\paragraph{Consistency and the signal-to-noise ratio} Under $H_0$
	($\Delta_i\equiv0$), the expectation and standard deviation are both of order
	$\sqrt N$ --- the latter by Theorem~\ref{thm:Var-TVC}. Under any fixed
	alternative $\mathbf p\ne\mathbf r$,
	the folded-normal expectation theorem \cite[Theorem~6.5]{banic_elezovic_p3}
	gives
	\[
		\E(d_{\mathbf r})=N\sum_i|\Delta_i|+O(\sqrt N)=\Theta(N).
	\]
	The separation measured in units of the null standard deviation is therefore
	\begin{equation}\label{eq:dXY-effect-size}
		\frac{\E_{H_1}(d_{\mathbf r}) - \E_{H_0}(d_{\mathbf r})}{\sqrt{\Var_{H_0}(d_{\mathbf r})}}
		= \frac{\Theta(N)}{\Theta(\sqrt N)} = \Theta(\sqrt N) \to \infty.
	\end{equation}
	The consistency claim itself follows directly from the multinomial law of
	large numbers~\cite{billingsley}:
	\[
		\frac{d_{\mathbf r}}N\ \longrightarrow\ \sum_i|\Delta_i|
		\quad\hbox{in probability}
	\]
	(cf.\ \cite[Corollaries~6.12--6.13]{banic_elezovic_p3} for the corresponding
	threshold-test statement).
	The limit is positive whenever $\mathbf p\ne\mathbf r$, because equal gradient
	profiles together with normalization would force $\mathbf p=\mathbf r$; under
	$H_0$ the limit is zero. Thus any null-calibrated threshold of order
	$\sqrt N$ is exceeded with probability tending to one. The separation grows
	as $\sqrt N$ with the signal $\sum_i|\Delta_i|$, the $\ell_1$ distance of the
	discrete derivatives.

	\paragraph{Numerical validation} The closed form (\ref{eq:Var-dXY-closed}) and
	the per-edge profile (\ref{eq:dXY-peredge-var}) were checked against exact and
	Monte Carlo variances; the results are collected in Table~\ref{tab:dXY-var-MC}.
	\begin{table}[ht]
	\centering
	\begin{tabular}{l|r|r|r|c}
		$(\mathbf p, \mathbf r)$ & $N$ & $\Var(d_{\mathbf r})$ & closed \eqref{eq:Var-dXY-closed} & rel.\ error \\\hline
		comb 9:1 vs uniform  & 1000 & $1290.000$  & $1290.000$ & $<10^{-6}\%$ \\
		                     & 8000 & $10316.8$   & $10320.0$  & $+0.03\%$ \\
		comb 9:1 vs comb 2:1 & 8000 & $10321.2$   & $10320.0$  & $-0.01\%$ \\
		geometric vs uniform & 1000 & $243.497$   & $241.6$    & $-0.78\%$ \\
		                     & 8000 & $1935.8$    & $1932.8$   & $-0.16\%$ \\
	\end{tabular}
	\caption{The closed form \eqref{eq:Var-dXY-closed} against $\Var(d_{\mathbf r})$, $n=8$. The $N=1000$ rows
	are exact: $1290.000$ for the comb rows (the linearization \eqref{eq:Var-dXY-closed} is exact up to an
	exponentially small term, so the closed form \emph{is} the variance there), and $243.497$ for the geometric
	row by a two-dimensional fast Fourier transform of the exact edge law, confirmed to the printed digits by
	$2\cdot10^7$ Monte Carlo replications. The $N=8000$ rows are Monte Carlo ($K=6\cdot10^5$; relative standard
	error $\sqrt{2/K}\approx0.2\%$). The comb has high$:$low mass ratio $9:1$ ($p_{\text{high}}=0.225$,
	$p_{\text{low}}=0.025$); the geometric sampling law has common ratio $\lambda=0.398$ ($p_1\approx0.60$), the
	reference being uniform.}
	\label{tab:dXY-var-MC}
	\end{table}
	The comb (all $|a_i|\ge8$) matches to the printed digits already at $N=1000$; the geometric
	law has small tail bins ($Np_8\approx1$ at $N=1000$) and three edges with $|a_i|<2.5$, so the
	exponentially small error of Theorem~\ref{thm:Var-dXY} is not yet in force there
	(Proposition~\ref{prop:dXY-mixed} evaluated at the actual $a_i$ gives $246.0$, $+1.0\%$, the central
	limit theorem itself being poor for a bin holding one count).

\section{Variance for small numbers of bins ($n\le4$)}
\label{sec:small-n}

	The theory above assumes $n\ge5$, and the restriction has a single source:
	the bivariate estimate of Proposition~\ref{prop:pair-O1} requires that the
	bins touched by a pair of edges do not exhaust the sample, $q_0>0$, because
	$q_0>0$ is what makes the two-dimensional increment walk aperiodic on the
	full lattice $\Z^2$. As recorded in Remark~\ref{rem:pair-boundary}, exactly
	two configurations violate this: the adjacent pair at $n=3$ and the disjoint
	pair at $n=4$. The case $n=2$ fails one step earlier --- its single edge has
	$p_1+p_2=1$, the increment has span two, and already the per-edge input
	\cite[Lemma~4.1]{banic_elezovic_p3} does not apply. This section settles all
	three cases.

	Each degeneration is visible in an elementary identity:
	\begin{equation}\label{eq:three-identities}
		d_{\mathbf p}\big|_{n=2}=2|X_1-Np_1|,\ \ \
		D_1-D_2\big|_{n=3}=3X_2-N,\ \ \
		D_1+D_3\big|_{n=4}=N-2(X_1+X_3).
	\end{equation}
	At $n=2$ the statistic reduces to a binomial absolute deviation; at $n=3$
	the pair $(D_1,D_2)$ is confined to a coset of an index-$3$ sublattice of
	$\Z^2$, and the kink $(N\delta_1,N\delta_2)$ sits at a position in the coset that moves
	with $\{Np_2\}$; at $n=4$ the disjoint pair $(D_1,D_3)$ is confined to a coset of the
	parity sublattice, with phase $\{N(p_1+p_3)\}$. This is the general-law counterpart of the uniform small-$n$ hierarchy
	of \cite[\S8]{banic_elezovic_p2}: each covariance block --- per-edge,
	adjacent, disjoint --- degenerates precisely at its minimal $n$, where the
	conditioning total is deterministic and no remainder bin is left to mix
	over. (The leading coefficient alone already follows for every $n\ge2$ from
	Corollary~\ref{cor:weak-limit}; the content of \S\ref{ssec:small-n34} is the $O(1)$
	remainder and the location of the oscillations.) Throughout,
	$b(k;N,p)=\binom Nkp^k(1-p)^{N-k}$ and $B(\ell;\kappa,\vartheta)$ is the
	binomial CDF with the conventions of
	\cite[Theorem~3.2]{banic_elezovic_p3}.

\subsection{The two-bin case}
\label{ssec:small-n2}

	At $n=2$ the centered statistic reduces to a folded binomial,
	$d_{\mathbf p}=2|X_1-Np|$ with $X_1\sim\Bin(N,p)$, and the reference-centered statistic to
	$d_{\mathbf r}=2|X_1-Nr|$: the fold sits at the mean in the first case and at the
	hypothesized mean in the second. The mean-absolute-deviation asymptotics of the first
	are developed in~\cite[Theorem~4.2]{elezovic_mad}, and the second is the subject
	of~\cite{elezovic_folded}, which also records the mean-centered variance; we record
	here only the formulas needed below.

	\emph{The centered statistic.} Here the second moment is exact, because
	$d_{\mathbf p}^2=\tilde D_1^2$ and hence
	$\E(d_{\mathbf p}^2)=4Npq=\Var(D_1)$; the variance therefore follows from De~Moivre's
	closed form for the mean (\cite{diaconis_zabell}; \cite[Proposition~4.1]{elezovic_mad}).
	Write $\nu=\lceil Np\rceil$ and
	$h_N=\nu-Np\in[0,1)$. Then, by \cite[Corollary~2.2 and Eq.~(2.6)]{elezovic_folded}
	(the expansion resting on \cite[Theorem~4.2]{elezovic_mad}),
	\begin{equation}\label{eq:n2-var}
	\begin{aligned}
		\Var(d_{\mathbf p})&=4Npq-16\,\nu^2q^2\,b(\nu;N,p)^2\\
		&=4pq\Bigl(1-\frac2\pi\Bigr)N-\frac{4pq}{3\pi}
		+\frac{8}{\pi}B_2(h_N)+O(N^{-1}),
	\end{aligned}
	\end{equation}
	$B_2(t)=t^2-t+\tfrac16$. The leading formula of Theorem~\ref{thm:Var-TVC}
	therefore continues formally to $n=2$, where there are no pairs:
	$\Var(d_{\mathbf p})=(1-\tfrac2\pi)\Var(D_1)+O(1)$. The two-bin formula gives the sharper
	information:
	the $O(1)$ becomes the explicit constant of \eqref{eq:n2-var}, and it
	\emph{oscillates}, through the lattice defect of the mean. This is the
	general-law form of the phenomenon that only $n=2$ lacks a single constant
	$\beta$ \cite[\S8.3]{banic_elezovic_p2}: at $p=\tfrac12$ the oscillation is the
	parity of $N$ and \eqref{eq:n2-var} returns
	$(1-\tfrac2\pi)N+(-1)^N/\pi+O(N^{-1})$, the two-bin uniform expansion of
	\cite[\S8.1]{banic_elezovic_p2}.

	\emph{The reference-centered statistic.} This is the two-bin case of
	Theorem~\ref{thm:Var-dXY}, deferred there. Exactly,
	$\Var(d_{\mathbf r})=4[Npq+N^2(p-r)^2]-(\E\,d_{\mathbf r})^2$; and for fixed $r\ne p$
	(after relabeling, $r<p$), by \cite[Proposition~3.7 and Eq.~(3.19)]{elezovic_folded},
	\begin{equation}\label{eq:n2-var-dXY}
		\Var(d_{\mathbf r})
		=4Npq-16\,N(p-r)\,b(m;N,p)\,S_0(\theta_N)\,\bigl(1+O(N^{-1})\bigr),
	\end{equation}
	with $q=1-p$, $\bar r=1-r$, $m=\lfloor Nr\rfloor$, $\theta_N=\{Nr\}$,
	$\eta=rq/(\bar rp)<1$ and
	$S_0(\theta)=\theta/(1-\eta)+\eta/(1-\eta)^2$. The leading term $4Npq=\Var(D_1)$
	is exactly the quadratic form of Theorem~\ref{thm:Var-dXY} evaluated at $n=2$
	(there $c_1=-\epsilon_1$, $c_2=\epsilon_1$), which proves the deferred case; and
	\eqref{eq:n2-var-dXY} sharpens its $O(e^{-c'N})$ remainder to an explicit
	constant, rate and oscillation, the rate being the Kullback--Leibler divergence
	$D(r\Vert p)$. The deficit is positive, as it must be: $|\cdot|$ is
	$1$-Lipschitz, so folding can only remove variance, and what it removes is the
	lower-tail crossing probability beyond the fold.

\subsection{The sublattice extension at $n=3$ and $n=4$}
\label{ssec:small-n34}

	We now treat the two pair configurations excluded by $q_0>0$. Throughout,
	$(j,k)$ is either the adjacent pair $(1,2)$ at $n=3$ or the disjoint pair
	$(1,3)$ at $n=4$, all $p_i>0$, and $\zeta$ denotes the per-trial increment of
	$(D_j,D_k)$.
	For the adjacent three-bin case set
	\[
		v_0=(0,1),\qquad
		\Lambda_3:=\{(x,y)\in\Z^2:\ x+2y\equiv0 \pmod 3\},
	\]
	and for the disjoint four-bin case set
	\[
		v_0=(0,1),\qquad
		\Lambda_2:=\{(x,y)\in\Z^2:\ x+y\equiv0\pmod2\}.
	\]

\begin{lemma}[Minimal lattices and cosets]\label{lem:lattices}
	\emph{(i)} For the adjacent pair $(D_1,D_2)$ at $n=3$,
	\[
		\operatorname{supp}\zeta\subset v_0+\Lambda_3,
		\qquad (D_1,D_2)\in Nv_0+\Lambda_3,
	\]
	$\Lambda_3$ is the minimal lattice generated by the differences of the support,
	and
	\[
		\Lambda_3^*=2\pi\Z^2\ \cup\
		\Bigl(\pm\tfrac{2\pi}3(1,-1)+2\pi\Z^2\Bigr),
	\]
	with $|\psi(\theta)|<1$ for $\theta\notin\Lambda_3^*$.

	\emph{(ii)} For the disjoint pair $(D_1,D_3)$ at $n=4$,
	\[
		\operatorname{supp}\zeta\subset v_0+\Lambda_2,
		\qquad (D_1,D_3)\in Nv_0+\Lambda_2,
	\]
	$\Lambda_2$ is the minimal lattice generated by the differences of the support,
	and
	$\Lambda_2^*=2\pi\Z^2\cup\bigl((\pi,\pi)+2\pi\Z^2\bigr)$, with
	$|\psi(\theta)|<1$ off $\Lambda_2^*$.

	\emph{(iii)} In both cases the support of $\zeta$ contains three affinely
	independent points of positive probability, so the per-trial covariance is
	nonsingular and $|\rho_{j,k}|<1$.
\end{lemma}

\noindent
Equivalently, the coset statements give the congruences
$D_1-D_2\equiv-N\pmod3$ at $n=3$ and
$D_1+D_3\equiv N\pmod2$ at $n=4$; these are the identities
\eqref{eq:three-identities} and $D_1+D_3=N-2(X_1+X_3)$ written as lattice
constraints.

\begin{proof}
	(i) The linear form $L(x,y)=x+2y$ takes the value $-1\equiv2$, $-1\equiv2$,
	$2$ on the three support points, so $\mathrm{supp}\,\zeta\subseteq
	v_0+\Lambda_3$ with $L(v_0)=2$. The differences
	$(2,-1)=(1,-1)-(-1,0)$ and $(1,1)=(0,1)-(-1,0)$ lie in $\Lambda_3$ and have
	determinant $3=[\Z^2:\Lambda_3]$, so they generate a sublattice of
	$\Lambda_3$ of the same covolume, hence all of $\Lambda_3$; the minimal
	lattice, which contains them and is contained in $\Lambda_3$, is therefore
	$\Lambda_3$. For the dual: $\frac{2\pi}3(1,-1)\cdot(x,y)=\frac{2\pi}3(x-y)$
	and $x-y\equiv x+2y\equiv0\pmod 3$ on $\Lambda_3$, so the displayed set lies
	in $\Lambda_3^*$; it has covering index $3$ over $2\pi\Z^2$, matching
	$\mathrm{covol}(\Lambda_3^*)=(2\pi)^2/\mathrm{covol}(\Lambda_3)$, so it is
	all of $\Lambda_3^*$. Finally
	$\psi(\theta)=p_1e^{-i\theta_1}+p_2e^{i(\theta_1-\theta_2)}+p_3e^{i\theta_2}$
	is a convex combination of unimodular numbers with positive weights, so
	$|\psi(\theta)|=1$ forces the three phases to coincide:
	$-\theta_1\equiv\theta_1-\theta_2\equiv\theta_2\pmod{2\pi}$, whence
	$3\theta_1\equiv0$ and $\theta_2\equiv-\theta_1$; these are exactly the
	points of $\Lambda_3^*$ modulo $2\pi\Z^2$.

	(ii) Identical scheme: $x+y$ is odd on all four support points;
	$(1,1),(1,-1)$ are differences with determinant $-2$; and
	$|\psi(\theta)|=1$ forces
	$-\theta_1\equiv\theta_1\equiv-\theta_2\equiv\theta_2$, i.e.\
	$\theta\in\{(0,0),(\pi,\pi)\}$ modulo $2\pi\Z^2$.

	(iii) At $n=3$ the three support points are affinely independent; at $n=4$,
	$(-1,0),(1,0),(0,1)$ are.
\end{proof}

\begin{lemma}[Fractional modes avoid the axes]\label{lem:frac-modes}
	Every $\xi\in\Lambda_3^*\setminus2\pi\Z^2$ has both coordinates nonzero,
	with $|\xi_1|,|\xi_2|\ge\frac{2\pi}3$; every
	$\xi\in\Lambda_2^*\setminus2\pi\Z^2$ has both coordinates odd multiples of
	$\pi$.
\end{lemma}

\begin{proof}
	For $\Lambda_3^*$: $\xi=2\pi(a\pm\tfrac13,\,b\mp\tfrac13)$ with $a,b\in\Z$,
	and $a\pm\tfrac13$, $b\mp\tfrac13$ are never integers, in particular never
	zero, and are at distance $\ge\tfrac13$ from $0$. For $\Lambda_2^*$:
	$\xi=(\pi(2a+1),\pi(2b+1))$.
\end{proof}

	Lemma~\ref{lem:frac-modes} identifies where the additional dual modes can
	contribute. In Step~4 of the
	proof of Proposition~\ref{prop:pair-O1} the $O(1)$ lattice
	corrections are carried by the \emph{axis} modes $\xi=(2\pi r_1,0)$,
	$(0,2\pi r_2)$, while the \emph{double} modes ($\xi_1\xi_2\ne0$) contribute
	only $O(N^{-1})$. The sublattice degeneration enriches the dual by new
	modes, but all of them belong to the double-mode family and hence contribute
	only at order $N^{-1}$.

\begin{theorem}[Extension of the bivariate estimate to the excluded configurations]
	\label{thm:sublattice}
	Let $(j,k)$ be the adjacent pair at $n=3$ or the disjoint pair at $n=4$,
	with all $p_i>0$ fixed. Then, as $N\to\infty$,
	\begin{equation}\label{eq:sublattice-SM}
		\E\bigl[\,|\tilde D_j|\,|\tilde D_k|\,\bigr]
		=\frac2\pi\,\sigma_j\sigma_k
		\bigl[\sqrt{1-\rho_{j,k}^2}+\rho_{j,k}\arcsin\rho_{j,k}\bigr]+O(1),
	\end{equation}
	and consequently
	$\Cov(|\tilde D_j|,|\tilde D_k|)
	=\frac2\pi\sigma_j\sigma_k\,\Phi(\rho_{j,k})+O(1)$.
\end{theorem}

\noindent
Within the $O(1)$, the integer axis modes contribute the same lattice-correction
family as in the generic case, while the fractional dual modes contribute only
$O(N^{-1})$. Their phases are governed jointly by $\{N\delta_j\},\{N\delta_k\}$ and by
the position of the mean $(N\delta_j,N\delta_k)$ relative to the coset $Nv_0+\Lambda$,
i.e.\ by $\{Np_2\}$ in the adjacent $n=3$ case and by $\{N(p_1+p_3)\}$ in the disjoint
$n=4$ case; at the uniform law these reduce to $N\bmod3$ and to the parity of $N$.

\begin{proof}
	We follow Steps~1--4 of the proof of
	Proposition~\ref{prop:pair-O1}, indicating what changes. Write
	$\Lambda$ for the minimal lattice of Lemma~\ref{lem:lattices} ($\Lambda_3$
	or $\Lambda_2$), $v_0$ for the support shift, and
	$d:=[\Z^2:\Lambda]\in\{3,2\}$, which is also the covolume of $\Lambda$.
	Fix an integer basis matrix $T$ of $\Lambda$: for $n=3$,
	$T=\begin{pmatrix}1&2\\1&-1\end{pmatrix}$ (columns $(1,1),(2,-1)$,
	$\det T=-3$); for $n=4$, $T=\begin{pmatrix}1&1\\1&-1\end{pmatrix}$
	($\det T=-2$).

	\emph{Step 1$'$ (weighted local expansion on the coset).} Define
	$\zeta'$ by $\zeta=v_0+T\zeta'$. For $n=3$ this gives
	$\zeta'\in\{(-1,0),(-1,1),(0,0)\}$ with probabilities $p_1,p_2,p_3$; for
	$n=4$, $\zeta'\in\{(-1,0),(0,1),(-1,1),(0,0)\}$ with probabilities
	$p_1,p_2,p_3,p_4$. In both cases the differences of the support of $\zeta'$
	contain $(1,0)$ and $(0,1)$, so $\zeta'$ is $\Z^2$-minimal and aperiodic,
	with nonsingular covariance and all moments finite. The weighted local
	Edgeworth expansion invoked in Step~1 of the proof of
	Proposition~\ref{prop:pair-O1}
	(Bhattacharya--Rao~\cite[Theorem~22.1]{bhattacharya_rao})
	therefore applies to $S_N'=\sum_t\zeta_t'$. Since
	$(D_j,D_k)=Nv_0+TS_N'$, transforming back (densities pick up the Jacobian
	$|\det T|^{-1}$, and the standardizations of $S_N'$ and of $(D_j,D_k)$
	differ by fixed $N$-free linear maps, so the polynomial weights are
	equivalent with $N$-free constants) gives, with $f_N$ the bivariate
	Edgeworth pseudo-density of $(D_j,D_k)$ built from its own cumulants,
	\[
		\sup_{m\in Nv_0+\Lambda}\ (1+|v_m|^6)\,
		\bigl|\Pr\bigl((D_j,D_k)=m\bigr)-d\,f_N(m)\bigr|=o(N^{-2}),
	\]
	while $\Pr\bigl((D_j,D_k)=m\bigr)=0$ off the coset. The covolume factor $d$
	is forced by normalization: the coset carries one lattice point per area
	$d$. The Riemann-sum argument of Step~1 now runs over a grid of covolume
	$d$, so $\sum_{m\in\text{coset}}g(v_m)=O(\sigma_j\sigma_k/d)$, and with
	$h(y):=|y_1-c_j|\,|y_2-c_k|\,f_N(y)$, $c_j=N\delta_j$, $c_k=N\delta_k$,
	\begin{equation}\label{eq:step1-coset}
		\E\bigl[\,|\tilde D_j|\,|\tilde D_k|\,\bigr]
		=d\!\!\sum_{m\in Nv_0+\Lambda}\!\! h(m)+o(1).
	\end{equation}

	\emph{Step 2$'$ (Poisson summation over the coset).} The Fourier decay
	$\hat h(\xi)=O\bigl((1+\xi_1^2)^{-1}(1+\xi_2^2)^{-1}\bigr)$ proved in
	Step~2 of the proof of Proposition~\ref{prop:pair-O1} is a property of $h$
	alone and is unchanged. Poisson summation for the lattice $\Lambda$ of
	covolume $d$ \cite[Ch.~VII]{stein_weiss}, applied at the shift $v_N:=Nv_0$,
	gives
	\[
		d\!\!\sum_{m\in v_N+\Lambda}\!\! h(m)
		=\sum_{\xi\in\Lambda^*}\hat h(\xi)\,e^{i\xi\cdot v_N},
	\]
	the series absolutely convergent by the same decay (the dual
	$\Lambda^*\supset2\pi\Z^2$ is a lattice of covolume $(2\pi)^2/d$, so the
	mode count per ball grows only by the factor $d$).

	\emph{Step 3$'$ (zero mode).} For $\xi=0$ the covolume factors have
	canceled and $\hat h(0)$ is the same integral as in the generic case: the
	Gaussian part is the Nabeya moment \eqref{eq:sheppard-mardia}, the
	$N^{-1/2}$ Edgeworth part vanishes by the parity of $P_1$ against the even
	weight $|v_1v_2|\phi_\rho$, and the $N^{-1}$ part is $O(1)$ --- Step~3 of
	the proof of Proposition~\ref{prop:pair-O1} verbatim, using
	$|\rho_{j,k}|<1$ from Lemma~\ref{lem:lattices}(iii).

	\emph{Step 4$'$ (nonzero modes).} Split $\Lambda^*\setminus\{0\}$ into the
	integer part $2\pi\Z^2\setminus\{0\}$ and the fractional part
	$\Lambda^*\setminus2\pi\Z^2$. On the integer part $e^{i\xi\cdot v_N}=1$, and
	the axis-mode and double-mode estimates (Claims~\ref{cl:axis}--\ref{cl:nonzero}) of the
	proof of Proposition~\ref{prop:pair-O1} apply unchanged: the axis
	modes sum to $O(1)$ (the two-dimensional lattice oscillation in
	$\{N\delta_j\},\{N\delta_k\}$), the integer double modes to $O(N^{-1})$. On
	the fractional part, Lemma~\ref{lem:frac-modes} puts every mode in the
	double-mode family with both frequencies bounded away from zero
	($|\xi_i|\ge2\pi/3$, resp.\ $\ge\pi$). This is the key point in lowering the order: an
	axis mode has one vanishing frequency and its $1/\xi_1^2$ envelope sums to $O(1)$,
	whereas a double mode is damped in \emph{both} frequencies and its peak height is the
	pseudo-density at the corner, which is $O(N^{-1})$. Explicitly, the tensor-kink estimate
	of Claim~\ref{cl:double} gives for each such mode the leading coefficient
	$4f_N(c_j,c_k)/(\xi_1^2\xi_2^2)$ with
	$f_N(c_j,c_k)=O\bigl((\sigma_j\sigma_k)^{-1}\bigr)=O(N^{-1})$. The remainder has
	the same summable tensor envelope as in Claim~\ref{cl:double}, because the
	fractional modes are bounded away from both axes. Hence, summing over the
	fractional dual modes,
	\[
		\Bigl|\sum_{\xi\in\Lambda^*\setminus2\pi\Z^2}\hat h(\xi)\,e^{i\xi\cdot v_N}\Bigr|
		\;\le\; \frac{C}{N}\sum_{\xi\in\Lambda^*\setminus2\pi\Z^2}
		\frac{1}{\xi_1^2\,\xi_2^2}+o(N^{-1})
		\;=\;O(N^{-1}),
	\]
	the series converging because both $|\xi_1|$ and $|\xi_2|$ are bounded below on the
	fractional part. Together with the factor $e^{-i\xi\cdot c}$ of the double-mode leading
	term, their phases are $e^{i\xi\cdot(v_N-c)}$; for the fractional generator this is
	$e^{\mp2\pi iNp_2}$ (from $\xi=\pm\frac{2\pi}3(1,-1)$, $v_0=(0,1)$, $c=(N\delta_1,N\delta_2)$),
	resp.\ $e^{2\pi iN(p_1+p_3)}$ (from $\xi=(\pi,\pi)$), which reduce to the residue
	$e^{\mp2\pi iN/3}$ and the parity $(-1)^N$ at the uniform law: the coset oscillations
	announced in the statement, demoted to order $N^{-1}$.

	Summing the four steps proves \eqref{eq:sublattice-SM}; subtracting
	$\E|\tilde D_j|\,\E|\tilde D_k|=\frac2\pi\sigma_j\sigma_k+O(1)$ (the
	per-edge expansion \cite[Lemma~4.1]{banic_elezovic_p3} applies at every edge
	for $n\ge3$, since each edge has $s_j<1$) gives the covariance form.
\end{proof}

\begin{corollary}[The variance theorem for all $n\ge3$]\label{cor:n34}
	Theorem~\ref{thm:Var-TVC} holds with the hypothesis $n\ge5$
	replaced by $n\ge3$: for fixed $\mathbf p$ with all $p_i>0$,
	\begin{equation}\label{eq:Var-TVC-final}
		\Var(d_{\mathbf p})=\Bigl(1-\frac2\pi\Bigr)\sum_{j=1}^{n-1}\Var(D_j)
		+\frac4\pi\sum_{1\le j<k\le n-1}\sigma_j\sigma_k\,\Phi(\rho_{j,k})
		+O(1).
	\end{equation}
\end{corollary}

\begin{proof}
	For $n\ge5$ this is Theorem~\ref{thm:Var-TVC}. For
	$n=3,4$ every pair is either covered by
	Proposition~\ref{prop:pair-O1} (the two adjacent pairs at $n=4$
	have $q_0=p_4>0$, resp.\ $q_0=p_1>0$) or by
	Theorem~\ref{thm:sublattice} (the adjacent pair at $n=3$, the disjoint pair
	at $n=4$), and the per-edge terms carry $O(1)$ remainders by
	\cite[Lemma~4.1]{banic_elezovic_p3}. At $n=2$ the sharper
	formula \eqref{eq:n2-var} replaces the display.
\end{proof}

\begin{remark}[The uniform reduction and the fourth-cumulant term]\label{rem:uniform-n34}
	At the uniform law the corollary reproduces the small-$n$ facts of
	\cite[\S8]{banic_elezovic_p2}: for $n=4$, $\rho_{1,3}=0$ and the disjoint
	$O(1)$ constant is $-1/(4\pi)$, obtained there from the exact skeleton
	$h(N)$; the numerics below recover it, with the parity oscillation decaying like
	$N^{-1}$, in accordance with Theorem~\ref{thm:sublattice} and with the
	observation of \cite[\S8.3]{banic_elezovic_p2}
	that at $n=3,4$ the residue oscillation enters at order $N^{-1}$.
	Likewise, the joint fourth-cumulant portion
	$\kappa_{2,2}/(2\pi\tilde s_j\tilde s_k\sqrt{1-\rho_{j,k}^2})$ of the
	disjoint $O(1)$ remainder, isolated in
	Proposition~\ref{prop:disj-dominant}, is a smooth Edgeworth
	computation with no lattice input, so it identifies the same portion of the
	constant at $n=4$; the numerical example below is consistent with this.
\end{remark}

\subsection{Exact finite-$N$ reductions and numerical validation}
\label{ssec:small-n-exact}

	The expansions above are validated against exact values. Both small-$n$
	configurations admit one-dimensional exact reductions with $O(N)$ total
	cost; they are also of independent interest, since they expose the
	degenerate geometry directly.

\begin{lemma}[Two-kink binomial expectation]\label{lem:two-kink}
	Let $Y\sim\Bin(\kappa,\vartheta)$ and let $\alpha\le\beta$ be real. Then
	\begin{equation}\label{eq:two-kink}
		\E\,|Y-\alpha|\,|Y-\beta|
		=\kappa\vartheta(1-\vartheta)+(\kappa\vartheta-\alpha)(\kappa\vartheta-\beta)
		+2\!\!\sum_{\alpha<y<\beta}\!\!(y-\alpha)(\beta-y)\Pr(Y=y),
	\end{equation}
	and the window sum is a linear combination of truncated moments
	$\sum_{y\le m}y^i\Pr(Y=y)$, $i=0,1,2$, each a binomial CDF:
	\[
		\sum_{y\le m}\Pr(Y=y)=B(m;\kappa,\vartheta),\quad
		\sum_{y\le m}y\Pr(Y=y)=\kappa\vartheta\,B(m-1;\kappa-1,\vartheta),
	\]
	\[
		\sum_{y\le m}y(y-1)\Pr(Y=y)
		=\kappa(\kappa-1)\vartheta^2\,B(m-2;\kappa-2,\vartheta).
	\]
\end{lemma}

\noindent
Hence \eqref{eq:two-kink} costs $O(1)$ binomial CDF evaluations.

\begin{proof}
	Pointwise, $|u-\alpha||u-\beta|=(u-\alpha)(u-\beta)
	+2(u-\alpha)(\beta-u)\mathbf 1_{(\alpha,\beta)}(u)$: outside $(\alpha,\beta)$
	the two factors have equal signs and the indicator term vanishes; inside,
	$(u-\alpha)(u-\beta)<0$ and adding twice its negative flips the sign. Taking
	expectations and using
	$\E(Y-\alpha)(Y-\beta)=\Var Y+(\E Y-\alpha)(\E Y-\beta)$ gives
	\eqref{eq:two-kink}. Expanding $(y-\alpha)(\beta-y)=-y^2+(\alpha+\beta)y
	-\alpha\beta$ reduces the window sum to truncated moments over an integer
	window, i.e.\ to differences of the displayed CDF identities, which follow
	from $y\binom\kappa y\vartheta^y=\kappa\vartheta\binom{\kappa-1}{y-1}
	\vartheta^{y-1}$ applied once and twice.
\end{proof}

\begin{proposition}[Exact adjacent mixed moment at $n=3$]\label{prop:n3-exact}
	Let $n=3$, all $p_i>0$, $\delta_1=p_2-p_1$, $\delta_2=p_3-p_2$,
	$\vartheta=p_1/(p_1+p_3)$. Conditionally on $X_2=m$, the variable
	$X_1$ is $\Bin(N-m,\vartheta)$ and \emph{both} centered edges are affine in
	it:
	\[
		\tilde D_1=-(X_1-\alpha_m),\qquad
		\tilde D_2=-(X_1-\beta_m),
	\]
	\[
		\alpha_m=m-N\delta_1,\qquad
		\beta_m=N-2m-N\delta_2,\qquad
		\alpha_m-\beta_m=3(m-Np_2).
	\]
	Let
	\[
		\alpha_m^-:=\min(\alpha_m,\beta_m),\qquad
		\alpha_m^+:=\max(\alpha_m,\beta_m).
	\]
	Consequently
	\begin{equation}\label{eq:n3-exact}
		\E\bigl[\,|\tilde D_1|\,|\tilde D_2|\,\bigr]
		=\sum_{m=0}^{N}b(m;N,p_2)\;
		\E\,\bigl|Y_m-\alpha_m\bigr|\,\bigl|Y_m-\beta_m\bigr|,
		\qquad Y_m\sim\Bin(N-m,\vartheta),
	\end{equation}
	where the inner expectation is evaluated by Lemma~\ref{lem:two-kink}
	with $(\alpha,\beta)=(\alpha_m^-,\alpha_m^+)$.
\end{proposition}

\noindent
Together with the exact per-edge mixture
\cite[Proposition~3.3]{banic_elezovic_p3}, this gives an exact evaluation of
$\Var(d_{\mathbf p})$ at $n=3$, requiring $O(N)$ binomial CDF evaluations.

\begin{proof}
	Given $X_2=m$, the remaining $N-m$ trials fall into bins $1$ and $3$ with
	renormalized probabilities $\vartheta$ and $1-\vartheta$
	(Proposition~\ref{prop:2bin-factor} applied to the pooled pair), so
	$X_1\sim\Bin(N-m,\vartheta)$ and $X_3=N-m-X_1$. Then
	$\tilde D_1=m-X_1-N\delta_1$ and $\tilde D_2=(N-m-X_1)-m-N\delta_2$, which
	are the
	displayed affine forms; the kink separation is
	$\alpha_m-\beta_m=3m-N(1+\delta_1-\delta_2)=3(m-Np_2)$, since
	$1+\delta_1-\delta_2=3p_2$ when $p_1+p_2+p_3=1$. Averaging over
	$X_2\sim\Bin(N,p_2)$ gives \eqref{eq:n3-exact}. Finally
	$\Var(|\tilde D_j|)=\sigma_j^2-(\E|\tilde D_j|)^2$ exactly, and
	$\Var(d_{\mathbf p})=\sum_j\Var(|\tilde D_j|)+2\Cov(|\tilde D_1|,|\tilde D_2|)$.
\end{proof}

\begin{remark}[The two kinks and the coset]\label{rem:two-kinks}
	Representation \eqref{eq:n3-exact} makes the index-$3$ geometry concrete:
	the pair $(\tilde D_1,\tilde D_2)$ is a \emph{deterministic} function of the
	single binomial $X_1$ once $X_2$ is fixed, and the two kinks sit at a
	separation $3(m-Np_2)$ that fluctuates on the scale $\sqrt N$ --- the
	one-dimensional shadow of the coset $D_1-D_2\equiv-N\pmod3$. This is the
	general-law analogue of the conditional trinomial skeleton used for the
	uniform adjacent pair in \cite[\S5.1]{banic_elezovic_p2}.
\end{remark}

\begin{proposition}[Exact disjoint covariance at $n=4$]\label{prop:n4-exact}
	Let $n=4$, all $p_i>0$. In the exact disjoint reduction
	Lemma~\ref{lem:disj-reduction} the mixing law degenerates: with
	$s_1=p_1+p_2$, $s_3=p_3+p_4=1-s_1$, the block totals satisfy $M_3=N-M_1$
	almost surely, so with $M:=M_1\sim\Bin(N,s_1)$,
	\begin{equation}\label{eq:n4-exact}
		\Cov\bigl(|\tilde D_1|,|\tilde D_3|\bigr)
		=\Cov\bigl(\hat G_1(M),\,\hat G_3(N-M)\bigr),
	\end{equation}
	where $\hat G_j(a)=\E\,|2\Bin(a,\pi_j)-a-N\delta_j|$,
	$\pi_j=p_{j+1}/s_j$, is the translated two-bin expectation, in the closed
	form of \cite[Theorem~3.2]{banic_elezovic_p3}.
\end{proposition}

\noindent
The right side of \eqref{eq:n4-exact} is a single sum over the binomial law of
$M$. Since each value of $\hat G_j$ costs $O(1)$ binomial CDF evaluations, the
total cost is $O(N)$.

\begin{proof}
	Immediate from Lemma~\ref{lem:disj-reduction}: the pair $(M_1,M_3)$ is
	trinomial with third-category mass $1-s_1-s_3=0$, hence supported on
	$M_1+M_3=N$. The conditional independence and the identification of the
	conditional folded means are as there; the closed form of $\hat G_j$ is
	\cite[Theorem~3.2]{banic_elezovic_p3}.
\end{proof}

	\paragraph{Numerical validation} All entries below are exact computations
	(no Monte Carlo): the $n=3$ values by Proposition~\ref{prop:n3-exact}
	(checked against full enumeration of the trinomial to $10^{-13}$), the
	$n=4$ values by Proposition~\ref{prop:n4-exact} (checked against full
	enumeration at $N=30$), the $n=2$ values by \eqref{eq:n2-var}.

	\emph{Two bins.} At $p=\tfrac12$ the exact $\Var(d_{\mathbf p})-(1-\frac2\pi)N$ equals
	$-0.3191$ and $+0.3175$ at $N=101,102$, against $\mp1/\pi=\mp0.3183$: the
	parity oscillation of \eqref{eq:n2-var}. For the reference-centered statistic at
	$p=0.5$, $r=0.42$, the deficit $4Npq-\Var(d_{\mathbf r})$ is $14.75$ at $N=100$ and
	$0.0084$ at $N=800$, and \eqref{eq:n2-var-dXY} reproduces the deficit to $14\%$ at
	$N=800$ and $1.8\%$ at $N=6400$, the relative error decaying like $\approx116/N$
	(relative to the variance it is $1.4\cdot10^{-6}$ at $N=800$).

	\emph{Three bins.} For $\mathbf p=(0.2,0.35,0.45)$ ($\rho_{1,2}=-0.5654$) the remainder
	$\Cov(|\tilde D_1|,|\tilde D_2|)-\tfrac2\pi\sigma_1\sigma_2\Phi(\rho_{1,2})$ equals $-0.0573$, $-0.0591$,
	$-0.0601$, $-0.0604$ at $N=60,120,240,360$ --- the values $N\equiv0\pmod{60}$, at which all lattice phases
	$\{N\delta_1\},\{N\delta_2\},\{Np_2\}$ are frozen, so along this subsequence the remainder converges like
	$c+O(N^{-1})$. Over consecutive $N$ it does not converge: for $N=240,\dots,261$ it ranges over
	$[-0.0886,-0.0601]$, a period-$20$ oscillation of peak-to-peak amplitude $0.029$, which is the $O(1)$
	axis-mode lattice term of Theorem~\ref{thm:sublattice}. The coset oscillation proper is visible where its
	phase $\{Np_2\}$ moves while the axis phases are frozen: at the uniform point
	$\mathbf p=(\tfrac13,\tfrac13,\tfrac13)$ ($\delta_j=0$, phase $N\bmod3$) the covariance remainder, after
	removal of its smooth part $c+a/N+b/N^2$ fitted on $N=240,\dots,252$, oscillates by $+8.6\cdot10^{-4}$
	($N\equiv0$) and $-5.4\cdot10^{-4}$ ($N\equiv1,2\pmod3$), the pattern $\cos(2\pi N/3)$ of order $N^{-1}$
	predicted by Step~4$'$; for $\mathbf p=(\tfrac29,\tfrac59,\tfrac29)$ and $N\equiv2\pmod3$ (axis phases frozen)
	the residual times $N$ has class means $+0.045,+0.047,-0.096$ over $N\bmod9$, confirming the phase
	$\{Np_2\}=\{5N/9\}$. For $\mathbf p=(0.3,0.3,0.4)$ and $N\in10\Z$ all three phases are frozen ($Np_2\in\Z$),
	and the remainder shows only the smooth drift ($-0.066284,-0.066305,-0.066325$ at $N=300,310,320$).

	\emph{Four bins.} At the uniform point the constant must be
	$-1/(4\pi)=-0.0795775$ \cite[\S8.2]{banic_elezovic_p2}; the exact values
	give
	\[
	\begin{array}{r|cccc}
		N & 100 & 101 & 400 & 401\\\hline
		\Cov(|\tilde D_1|,|\tilde D_3|)+\tfrac1{4\pi}
		& +3.0\cdot10^{-4} & -1.3\cdot10^{-3}
		& +7.5\cdot10^{-5} & -3.2\cdot10^{-4}
	\end{array}
	\]
	--- value confirmed, the parity oscillation decaying like $N^{-1}$ (a factor $4$ from $N=100$ to $400$;
	Step~4$'$ gives $(-1)^N/(4\pi N)$). For the non-uniform $\mathbf p=(0.1,0.3,0.25,0.35)$ ($\rho_{1,3}=-0.0434$)
	the remainder over the Nabeya term is $-0.0695$, $-0.0697$, $-0.0698$ at the phase-aligned $N=100,200,400$
	($N\equiv0\pmod{20}$); over $N=400,\dots,420$ it lies in $[-0.0702,-0.0698]$ (range $4.3\cdot10^{-4}$),
	$|R(401)-R(400)|=3.6\cdot10^{-4}$, $|R(801)-R(800)|=2.2\cdot10^{-4}$ --- not halving; the genuine parity
	component (even mean minus odd mean) is $2.3\cdot10^{-5}$ ($N\approx400$), $1.7\cdot10^{-5}$ ($N\approx800$),
	and the dominant part is the $O(1)$ axis-lattice term with phases $\{0.2N\},\{0.1N\},\{0.35N\}$. The
	fourth-cumulant portion of Proposition~\ref{prop:disj-dominant} equals $-0.0644$ here, consistent with its
	role as one identified part of the constant.

\section{Conclusion}

	The variance theory for the matched statistic $d_{\mathbf p}$ has the same folded-normal
	core as the uniform theory, but with an important new feature: disjoint edge
	pairs no longer disappear at leading order. The main formula expresses
	$\Var(d_{\mathbf p})$ through the exact multinomial edge variances and correlations and
	the Nabeya function $\Phi(\rho)$, with an $O(1)$ per-pair remainder
	at fixed positive probability vector. This recovers the uniform leading
	variance when $p_i\equiv1/n$, while showing why non-uniform reference laws with
	non-small gradients, comb-shaped ones in particular, require the disjoint terms in the
	null variance (under a comb \emph{alternative} to a different null the absolute values
	linearize instead, Theorem~\ref{thm:Var-dXY}).

	The profile calculations give two useful checks on the formula. Symmetric combs
	admit a closed expression in terms of two adjacent correlations and one
	disjoint magnitude, whereas positive $C^1$ profiles have, in the large-bin
	limit, the same universal constant as the uniform law. The
	reference-centered statistic $d_{\mathbf r}$ shows that fixed nonzero mismatches
	linearize the absolute values and reduce the leading variance to a multinomial
	quadratic form. The bivariate folded-normal covariance then
	interpolates between the matched $d_{\mathbf p}$ regime and this fully linearized
	alternative regime.

	The closing section removes the small-$n$ restriction: the main
	variance formula holds for every $n\ge3$, and at $n=2$ in a sharper, fully
	explicit form. The small-$n$ hierarchy mirrors the uniform
	case~\cite{banic_elezovic_p2}: each covariance block keeps, at its minimal
	$n$, a residue of the lattice on which its conditioning degenerates, but
	under a general law, exactly as under the uniform one, only the per-edge
	block at $n=2$ retains an additional oscillation at constant order. At
	$n=3$ and $n=4$ the generic integer-axis lattice corrections remain inside
	the fixed-\(\mathbf p\) \(O(1)\) remainder, while the extra coset oscillations from the
	fractional dual modes (with phases $\{Np_2\}$ and $\{N(p_1+p_3)\}$, which are $N\bmod3$
	and the parity of $N$ at the uniform law) are pushed to order $N^{-1}$: those modes have
	both coordinates nonzero.

\section*{Acknowledgments}
	During the preparation of this manuscript, the author used an AI assistant
	(Claude Opus 5.5, Anthropic) for computational and proof verification, 
	literature searches, and language editing. The author wrote the manuscript, 
	reviewed all AI-assisted output, and takes full responsibility for its content.

	\paragraph{Code and data} The scripts behind the tables and the exact small-$n$
	computations will be made public together with a monograph in preparation; until
	then they are available from the authors on request.



\begin{thebibliography}{99}

\bibitem{hoel1946}
	P.~G. Hoel,
	``The efficiency of the mean moving range,''
	\emph{Ann. Math. Statist.}, vol.~17, no.~4, pp.~475--482, 1946, doi:10.1214/aoms/1177730886.

\bibitem{kamat1953}
	A.~R. Kamat,
	``On the mean successive difference and its ratio to the root mean square,''
	\emph{Biometrika}, vol.~40, no.~1/2, pp.~116--127, 1953, doi:10.1093/biomet/40.1-2.116.

\bibitem{diaconis_zabell}
	P.~Diaconis and S.~Zabell,
	``Closed Form Summation for Classical Distributions: Variations on a Theme of
	De Moivre,''
	\emph{Statist. Sci.}, vol.~6, no.~3, pp.~284--302, 1991.

\bibitem{banic_elezovic_p5}
	N.~Bani\'c and N.~Elezovi\'c,
	``Discrete total variation under the multinomial, Part~5: cyclic total variation,'' in preparation, 2026.

\bibitem{banic_elezovic_tvor}
	N.~Bani\'c and N.~Elezovi\'c,
	``TVOR: Finding Discrete Total Variation Outliers Among Histograms,''
	\emph{IEEE Access}, vol.~9, pp.~1807--1832, 2021, doi:10.1109/ACCESS.2020.3047342.

\bibitem{banic_elezovic_ct}
	N.~Bani\'c and N.~Elezovi\'c,
	``Comb Test: Histogram Uniformity Testing Based on Discrete Total Variation,''
	\emph{IEEE Signal Process. Lett.}, vol.~33, pp.~3726--3730, 2026,
	doi:10.1109/LSP.2026.3728406 (arXiv:2606.01465).

\bibitem{vonneumann1941}
	J.~von Neumann,
	``Distribution of the Ratio of the Mean Square Successive Difference to the Variance,''
	\emph{Ann. Math. Statist.}, vol.~12, no.~4, pp.~367--395, 1941.

\bibitem{durbin_watson1950}
	J.~Durbin and G.~S. Watson,
	``Testing for Serial Correlation in Least Squares Regression. I,''
	\emph{Biometrika}, vol.~37, no.~3/4, pp.~409--428, 1950.

\bibitem{banic_elezovic_p1}
	N.~Bani\'c and N.~Elezovi\'c,
	``Discrete total variation under the multinomial, Part~1: expectation under the
	uniform null,'' preprint, \texttt{arXiv:2610.13366}, 2026.

\bibitem{banic_elezovic_p2}
	N.~Bani\'c and N.~Elezovi\'c,
	``Discrete total variation under the multinomial, Part~2: variance under the
	uniform null,'' preprint, \texttt{arXiv:2610.17046}, 2026.

\bibitem{elezovic_mad}
	N.~Elezovi\'c,
	``Local binomial expansions with an Appell shift, and the mean absolute deviation
	of the binomial distribution,'' Acta Math. Hungar., to appear; \texttt{arXiv:2607.18494}.

\bibitem{elezovic_folded}
	N.~Elezovi\'c,
	``Absolute deviations of the binomial about a prescribed centre: the tail expansion
	in closed form,'' preprint, \texttt{arXiv:2609.19064} [math.PR], 2026.

\bibitem{banic_elezovic_p3}
	N.~Bani\'c and N.~Elezovi\'c,
	``Discrete total variation under the multinomial, Part~3: expectation under a
	general distribution,'' preprint, \texttt{arXiv:2610.41255}, 2026.

\bibitem{bhattacharya_rao}
	R.~N. Bhattacharya and R.~Ranga~Rao,
	\emph{Normal Approximation and Asymptotic Expansions}.
	Wiley, New York, 1976.

\bibitem{billingsley}
	P.~Billingsley,
	\emph{Probability and Measure}, 3rd ed.
	Wiley, New York, 1995.

\bibitem{johnson_kotz_balakrishnan}
	N.~L. Johnson, S.~Kotz, and N.~Balakrishnan,
	\emph{Discrete Multivariate Distributions}.
	Wiley, New York, 1997.

\bibitem{kanrobotti2017}
	R.~Kan and C.~Robotti,
	``On Moments of Folded and Truncated Multivariate Normal Distributions,''
	\emph{J. Comput. Graph. Statist.}, vol.~26, no.~4, pp.~930--934, 2017, doi:10.1080/10618600.2017.1322092.

\bibitem{leone1961}
	F.~C. Leone, L.~S. Nelson, and R.~B. Nottingham,
	``The Folded Normal Distribution,''
	\emph{Technometrics}, vol.~3, no.~4, pp.~543--550, 1961.

\bibitem{nabeya1951}
	S.~Nabeya,
	``Absolute moments in 2-dimensional normal distribution,''
	\emph{Ann. Inst. Statist. Math.}, vol.~3, no.~1, pp.~2--6, 1951.

\bibitem{pearson1900}
	K.~Pearson,
	``On the Criterion that a Given System of Deviations from the Probable in the Case
	of a Correlated System of Variables is Such that it Can Be Reasonably Supposed to Have
	Arisen from Random Sampling,''
	\emph{Philos. Mag. (5)}, vol.~50, no.~302, pp.~157--175, 1900.

\bibitem{patnaik1949}
	P.~B. Patnaik,
	``The Non-Central $\chi^2$- and $F$-Distributions and their Applications,''
	\emph{Biometrika}, vol.~36, no.~1/2, pp.~202--232, 1949.

\bibitem{price1958}
	R.~Price,
	``A Useful Theorem for Nonlinear Devices Having Gaussian Inputs,''
	\emph{IRE Trans. Inform. Theory}, vol.~4, no.~2, pp.~69--72, 1958, doi:10.1109/TIT.1958.1057444.

\bibitem{sheppard1899}
	W.~F. Sheppard,
	``On the Application of the Theory of Error to Cases of Normal Distribution and
	Normal Correlation,''
	\emph{Philosophical Transactions of the Royal Society of London A}, vol.~192, pp.~101--167, 1899.

\bibitem{stein_weiss}
	E.~M. Stein and G.~Weiss,
	\emph{Introduction to Fourier Analysis on Euclidean Spaces}.
	Princeton University Press, Princeton, 1971.

\end{thebibliography}
\end{document}